\documentclass{article}

\usepackage[a4paper, total={5.5in, 8in}]{geometry}
\usepackage[T1]{fontenc}
\usepackage[utf8]{inputenc}
\usepackage{comment}
\usepackage{amsfonts}
\usepackage{amsmath}
\usepackage{amsthm}
\usepackage{pgfplots}
\usepackage{tikz}
\usepackage{mathabx}
\usepackage{comment}
\usepackage{graphicx,subcaption}
\usepackage{wrapfig}

\usetikzlibrary{positioning} 
\usetikzlibrary{calc}
\usetikzlibrary{arrows,shapes,backgrounds}

\newtheorem{theorem}{Theorem}[section]
\newtheorem{definition}[theorem]{Definition}
\newtheorem{remark}[theorem]{Remark}
\newtheorem{lemma}[theorem]{Lemma}

\newtheorem{proposition}[theorem]{Proposition}

\begin{document}
	
\title{Ill-posedness of a quasilinear wave equation in two dimensions for data in $H^{\frac{7}{4}}$}

\author{Gaspard Ohlmann \\ University of Basel, \emph{gaspard.ohlmann@unibas.ch}}
\maketitle

\begin{center}
	\textbf{Abstract}
\end{center}

\emph{In this article, we study the ill-posedness of a quasilinear wave equation. It was shown by Tataru and Smith in 2005 that for any $s>7/4$ (or $11/4$ in our situation), the equation is well-posed in $H^{s}\times H^{s-1}$. We show a sharpness result by exhibiting a quasilinear wave equation and an initial data such that the Cauchy problem is ill-posed for in $H^{11/4} (\ln H)^{-\beta}\times H^{7/4} (\ln H)^{-\beta}$. }
	
	\section{Introduction}

We study the well-posedness of quasilinear wave equations. We will consider the following equation

\begin{equation}\label{QLW}
	\sum_{i,j=0}^n g^{ij}(u,u') \partial_{x^i} \partial_{x^j} u = F(u,u'), \hspace{0.2cm} (t,x) \in S_T = [0,T[\times \mathbb{R}^n,
\end{equation}
where $\partial_{x^0} = \partial_t$ and $G=(g^{ij})$ and $F$ are smooth functions. Also, we assume that $g$ is close to the Minkowski metric $m$ ; i.e.,

\begin{equation}\label{elli}
	\sum_{i,j=0}^n |g^{ij} - m^{ij}| \leq 1/2.
\end{equation}

We will also define the corresponding Cauchy problem, 

\begin{equation}\label{QLW1.5}
	\left\{
	\begin{split}
		&\sum_{i,j=0}^n g^{ij}(u,u') \partial_{x^i} \partial_{x^j} u = F(u,u'), \hspace{0.2cm} (t,x) \in S_T = [0,T[\times \mathbb{R}^n,\\
		&(u,\partial_t u)_{|t=0} = (f,h),
	\end{split}
	\right.
\end{equation}
where $\partial_{x^0} = \partial_t$ and $G=(g^{ij})$ and $F$ are smooth functions.

Informally, the concept of well-posedness usually involves the existence, uniqueness, and continuity with respect to the initial condition. The weakest assumption we have to make on the initial data for the quasilinear equation to be well-posed is a subject that is widely studied for its mathematical interest as well as its application to physics. It is common to assume that the initial condition, as well as its derivative respectively belong to some Sobolev spaces $H^s$ and $H^{s-1}$. 

Using Sobolev estimate, one can show that for $s\geq n+1$, the problem is well-posed. 

An improvement has for instance been made by H. Bahouri and J. Y. Chemin in \cite{bahouri1999equations}, and in parallel by D. Tataru in \cite{tataru2000strichartz}, where they show the well-posedness of the equation for $s \geq \frac{n+1}{2} + \frac{1}{4}$, using mainly Strichartz estimates. 

Further improvements have been made, and the best result for well-posedness in low dimension has been made by D. Tataru and H. F. Smith in 2005 in \cite{smith2005sharp} and is described by the next theorem. Note that here, the metric is allowed to depend on $u$ but not on its derivatives.

\begin{theorem}\cite{smith2005sharp}\label{smithsharp}
	
	We consider the following Cauchy problem. 
	
	\begin{equation}\label{QLW2}
		\left\{
		\begin{split}
			&\sum_{i,j} g^{ij}(u) \partial_{x^i} \partial_{x^j} u = \sum_{i,j} q^{ij} (u) \partial_{x^i} u \partial_{x^j} u, \hspace{0.2cm} (t,x) \in S_T = [0,T[\times \mathbb{R}^n,\\
			&(u,\partial_t u)_{|t=0} = (f,g),
		\end{split}
		\right.
	\end{equation}
	where $\partial_{x^0} = \partial_t$ and $G=(g^{ij})$ and $Q=(q^{ij})$ are smooth functions. Also, we assume that $g$ is close to the Minksowski metric $m$.
	
	The Cauchy problem (\ref{QLW2}) is locally well-posed in $H^s \times H^{s-1}$ provided that 
	\begin{equation}
		\begin{split}
			s &> \frac{n}{2} + \frac{3}{4} =\frac{7}{4} \hspace{0.3cm} \text{ for } n=2, \\
			s &> \frac{n+1}{2} \hspace{0.35cm} \text{ for } n=3,4,5. 
		\end{split}
	\end{equation}
\end{theorem}

In 2005, D. Tataru and P. Smith proved in \cite{smith2005sharp} that the quasilinear wave equation (\ref{QLW1.5}) where $g$ is not allowed to depend on the gradient of $u$, is well-posed provided that $s>7/4$. In this paper, we show that we can find an initial condition belonging to $H^{11/4}(\ln H)^{-\beta} \times H^{7/4} (\ln H)^{-\beta}$ such that the equation (\ref{QLW1.5}) is ill-posed in dimension $2+1$, in the sense that there exists no time $T>0$ such that there exists a solution $u$ on $[0,T[$ with the property that $u(t,\cdot), \partial_t u(t,\cdot) \in H^{11/4} (\ln H)^{-\beta} \times H^{7/4} (\ln H)^{-\beta}$ for a $\beta>1/2$. We show later that this implies the ill-posedness of the equation considered in \cite{smith2005sharp} for initial conditions belonging to $u(t,\cdot), \partial_t u(t,\cdot) \in H^{7/4} (\ln H)^{-\beta} \times H^{3/4} (\ln H)^{-\beta}$ for a $\beta>1/2$. We will show that the logarithmically modified Sobolev space $H^{s} (\ln H)^{-\beta}$ contains function slightly less regular than $H^s$ when $\beta>0$, in the sense that $H^{s}\subseteq H^s (\ln H)^{-\beta} \subseteq H^{s-\lambda}$ for any $\lambda >0$.
Note that in fact, the initial condition that we create has a compact support. 

We aim to provide a negative result. This is done by finding an equation and an initial condition such that the resulting Cauchy problem is not well-posed, meaning that there exists no time $T$ such that the problem is well-posed on $[0,T[$. Such a phenomenon is known as instantaneous blow-up. A negative result for $n=3$ has already been done by H. Lindblad in \cite{lindblad1998counterexamples}. The corresponding index $s$ for his problem is $3$, he hence exhibits the sharpness of the criteria established in \cite{smith2005sharp}.

First, we need to compute the corresponding index when the functions $G$ and $Q$ are allowed to also depend on $\nabla u$ and not only on $u$.

\begin{lemma}\label{corindice}
	Introducing the Cauchy problem with a $\nabla u$ dependency on $G$ and $Q$ as the following
	
	\begin{equation}\label{QLW7}
		\left\{
		\begin{split}
			&\sum_{i,j} g^{ij}(u,\nabla u) \partial_{x^i} \partial_{x^j} u = \sum_{i,j} q^{ij} (u,\nabla u) \partial_{x^i} u \partial_{x^j} u, \hspace{0.2cm} (t,x) \in S_T = [0,T[\times \mathbb{R}^3, \\
			&(u,\partial_t u)_{|t=0} = (f,g),
		\end{split}
		\right.
	\end{equation}
	where $\partial_{x^0} = \partial_t$ and $G=(g^{ij})$ and $Q=(q^{ij})$ are smooth functions. Also, we assume that $g$ is close to the Minkowski metric $m$ ;
	
	The Cauchy problem (\ref{QLW7}) is locally well-posed in $H^s \times H^{s-1}$ provided that 
	\begin{equation}
		\begin{split}
			s &> \frac{n}{2} + \frac{3}{4} + 1 \hspace{0.3cm} \text{ for } n=2, \\
			s &> \frac{n+1}{2} + 1 \hspace{0.35cm} \text{ for } n=3,4,5. 
		\end{split}
	\end{equation}
	
\end{lemma}

\begin{proof} (Of lemma \ref{corindice})
	
	We consider the problem given by (\ref{QLW7}). Differentiating the equation with respect to $x_k$, we obtain
	
	\begin{multline}\label{systemlol}
		\sum_{i,j} \left[ \left( \partial_k u \partial_1 g^{ij}(u,\nabla u) \right) \partial_i \partial_j u + \sum_\alpha \left( \partial_{\alpha+1} g^{ij}(u,\nabla u) \right) \partial_k \partial_\alpha u  \cdot \partial_i \partial_j u  + \left( g^{ij}(u,\nabla u) \right) \partial_i \partial_j \partial_k u \right] \\
		= \sum_{i,j} \Bigg[ \left( \partial_k u \partial_1 g^{ij} (u,\nabla u) \right) \partial_i u \partial_j u + \sum_\alpha  \left( \partial_i u \partial_{\alpha+1} q^{ij}(u,\nabla u)\right) \partial_j u \cdot \partial_k \partial_\alpha u \\
		+ q^{ij} (u,\nabla u) \partial_i \partial_k u \cdot \partial_j u + q^{ij} (u,\nabla u) \partial_i  u \cdot \partial_j \partial_k u \Bigg]
	\end{multline}
	
	Consequently, the system (\ref{systemlol}) for every $k$ can be put in the form
	
	\begin{equation}\label{s-1}
		\sum_{i,j} \tilde{g}^{ij}_k (v) \partial_i \partial_j v = \sum_{ij} \tilde{q}^{ij}_k(v) \partial_i v \partial_j v,
	\end{equation}
	
	where $v = (u,\nabla u)$. Now, $v \in H^{s-1}$ so (\ref{s-1}) is well-posed.
	
\end{proof}

Two mechanisms can create a blow-up. One of them is a space-independent blow-up. In a nutshell, the blow-up is caused by the underlying ODE that itself leads to a solution that blows up. Another kind of blow-up is caused by the focusing of the characteristics to a single point, leading to an infinite increase of the derivatives. Informally, in the case of a blow-up at $t=T$, if we denote by $\phi$ a characteristic function and $u$ the solution of the equation, in the first case, we see a phenomenon of the form (up to derivatives)

\begin{equation}
	u(\phi) \xrightarrow[t\rightarrow T]{} \infty, \hspace{0.2cm} \phi \neq 0,
\end{equation}

happening, whereas in the second case, one observes a phenomenon of the form

\begin{equation}
	\partial \phi(\nu) \xrightarrow[t\rightarrow T]{} 0,
\end{equation}

for a certain $\nu$. This will also lead to a blow-up of the derivatives using the chain rule. Those types of blow-up have already been described by Alinhac for instance in \cite{alinhac1995blowup}.

We create a counterexample for the dimension $1+2$, keeping in mind the intuition provided in \cite{lindblad1998counterexamples}. Our goal is to create a second type (geometric) blow-up. The differences here, are that we want a function that is globally defined because since the characteristic function is not in the desired Sobolev space, it is not clear that we can extend a counterexample on a set $\Omega\subset \mathbb{R}^2$ to a counterexample that contradicts the well-posedness on $\mathbb{R}^2$ with the sharp index $s=\frac{11}{4}$ if $g$ is allowed to depend on $\nabla u$, or $s=\frac{7}{4}$ if $g$ only depends on $u$. To have a sound argument, we first consider a regularized initial condition that leads to a blow-up at a time that is getting smaller as the initial condition is being less regularized, and we then construct an initial condition that leads to an instantaneous blow-up. An additional difference is that the corresponding index here is not an integer, and this leads to many technical issues. Lastly, we introduce the logarithmic Sobolev spaces $H^{11/4}(\ln H)^{-\beta}$ and $H^{7/4}(\ln H)^{-\beta}$, because they are the set in which our data perfectly fit.

We consider the model equation and the corresponding Cauchy problem

\begin{equation}\label{modelequation}
	\left\{
	\begin{split}
		&\Box u = (D u) D^{2} u,  \\
		&(u,\partial_t u)_{|t=0} = (f,g),
	\end{split}
	\right.
\end{equation}
where $D= (\partial_{x_1} - \partial_{t})$.

Note that (\ref{modelequation}) is of the form (\ref{QLW}) with

\begin{equation}\label{metricg}
	g=
	\begin{bmatrix}
		1-v & v & 0  \\
		v & -1-v & 0  \\
		0 & 0 & -1  \\
	\end{bmatrix}, v= Du.
\end{equation}

\begin{remark}
	Using a scaling argument, we can see that if $s<2$, then the problem is ill-posed.
	
\end{remark}

The counterexample we produce is in a slightly less regular space than $H^{11/4}$. Indeed we will consider the logarithmic perturbation of $H^{11/4}$, denoted by $H^{11/4}(\ln H)^{-\beta}$ as the set of functions $f$ such that the $L^2$ norm of $|\xi|^{11/4} \left( 1+ |\ln(|\xi|)|^{-\beta} \right) \mathcal{F}(f)(\xi)$, where $\mathcal{F}$ denotes the Fourier transform, is finite. Our counterexample will belong to the set $H^{11/4}(\ln H)^{-\beta}$ with $\beta>1/2$. For the functions that we consider, we will show that this set is located between $H^{11/4-\varepsilon}$ and $H^{11/4}$, for any $\varepsilon>0$. Hence, this proves the optimality of the index $11/4$ in the context of usual Sobolev spaces. It is however interesting to notice that the function that we create can be defined as a function in $H^{11/4} \times H^{7/4}$, as it is done in Appendix \ref{appendixA}, for which we expect to witness the same behavior, but the method we use to show the blow-up can not be applied anymore. The argument that does not hold in this situation, and the reason why we need slightly less regularity for the proof to hold, lies in the fact that the point where the blow-up occurs can not be proven to be in the domain of dependence anymore.

In the following, we construct a solution (\ref{modelequation}) with initial data in $H^{11/4}(\ln H)^{-\beta}\times H^{7/4}(\ln H)^{-\beta}$, that blows up instantly at $t=0^+$, as formulated in the following theorem.

\begin{theorem}\label{maingeneral} 
	We consider $\beta>1/2$. There exists initial data $(f,g) \in \dot H^{11/4}(\ln H)^{-\beta} \times \dot H^{7/4}(\ln H)^{-\beta}$ supported on a compact set, with $||f||_{\dot H^{11/4}(\ln H)^{-\beta}} + ||g||_{\dot H^{7/4} (\ln H)^{-\beta}}$ arbitrarily small, such that (\ref{modelequation}) does not have any proper solution $u$ such that 
	
	$(u,\partial_t u) \in C\left( [0,T[; \dot H^{11/4}(\ln H)^{-\beta} (\mathbb{R}^2) \times \dot H^{7/4} (\ln H)^{-\beta}(\mathbb{R}^2) \right)$ for any $T>0$.
	
	In fact, for all $\lambda>0$ small enough, (\ref{modelequation}) does not have any proper solution $u$ such that \\ $(u,\partial_t u) \in C\left( [0,T[; \dot H^{11/4-\lambda} (\mathbb{R}^2) \times \dot H^{7/4-\lambda}(\mathbb{R}^2) \right)$ for any $T>0$.

\end{theorem}

\subsection*{Notations}

First, for an integrable function $f:\mathbb{R}^n \to \mathbb{R}$, we denote by $\hat{f}:\mathbb{R}^n \to \mathbb{R}$ the function 

\begin{equation}
	\mathcal{F}(f)(\xi) = \hat{f}(\xi) = \int_{\mathbb{R}^n} f(y) e^{ - 2 i \pi <y\cdot \xi>} dy.
\end{equation}

We also define the Sobolev norm denoted by $||\cdot|| _{H(\mathbb{R}^n)^s}$ as 

\begin{equation}\label{sobosobo}
	||f(\cdot)||_{\dot H^{s}(\mathbb{R}^{n})}^2 = \int_{\mathbb{R}^n} |\xi|^{2s} |\hat f (\xi)|^{2} d\xi,
\end{equation}

and the corresponding Sobolev space $\dot H^{s} (\mathbb{R}^n)$ of functions such that this norm is finite.

When $s$ is an integer, the following holds 

\begin{equation}
	(2\pi)^s \int_{\mathbb{R}^n} |\xi|^{2s} |\hat f (t,\xi)|^{2} d\xi = \int_{\mathbb{R}^n} \left| \nabla^s f(t,x) \right|^2 dx,
\end{equation}

where $\nabla^s u = (\partial^s_1 u, ... , \partial^s_n u). $

We also define the notion of domain of dependence and the corresponding notations.

\begin{definition}\label{dependence}
	Let $\Omega \subset \mathbb{R}_+ \times \mathbb{R}^2$
	be an open set equipped with a Lorentzian metric $g_{j,k}$ satisfying (\ref{elli}). It is a domain of dependence for $g$ if the closure of the causal past $\Lambda_{t',x'}$ of each $(t',x') \in \Omega$ is contained in $\Omega$, with $z\in \Lambda_{t',x'}$ iff it can be joined to $(t',x')$ by a Lipschitz continuous curve $(t,x(t))$ satisfying 
	\begin{equation}\label{condlight}
		\sum_{i,j=0}^2 g_{i,j} (x) \frac{dx_i}{dt} \frac{dx_j}{dt} \geq 0,
	\end{equation}
	almost everywhere.
\end{definition}

For a domain $\Omega \subseteq \mathbb{R}^{1+n}$, we denote by $\Omega_t$ the set
\begin{equation}
	\Omega_t = \{ (\tau,x) \in \Omega, \hspace{0.2cm} t=\tau  \}.
\end{equation}

Also, for a set $\Omega \subseteq \mathbb{R}^{1+2}$, we denote (the dependence in $\Omega$ is not explicitly written) 

\begin{equation}
	a_t(x) = \left| \{ y \in \mathbb{R}, \hspace{0.2cm} (t,x_1,x_2)\in \Omega  \} \right|.
\end{equation}

\begin{remark}
	If we consider the Lorentzian metric that defines the linear wave equation, the causal past of a point is its associated light cone.
\end{remark}

We will make free use of the Huygens principle, meaning that for a solution $u$ defined on a domain of dependence $\Omega$ corresponding to the Lorentzian metric involved in the equation, the values of $u$ on the set $\Omega_t$ only depend of the value of $u$ on the set $\Omega_0$. This will be useful as we will first solve the equation for a locally defined function, and later create an extension of the initial condition. The computation previously made will remain valid in the corresponding domain of dependence.

For $\alpha$ a multi-index, we will use the following definition for fractional derivative 

\begin{equation}\label{coming2}
	\frac{\partial^\alpha f}{\partial x^\alpha} (x) = \mathcal{F}^{-1} \left( (2 \pi i)^s \xi_x^\alpha \cdot \mathcal{F}( f)(\xi) \right)(x).
\end{equation}

\section{Strategy and control of the initial condition}\label{controlcontrol}

In this chapter, we quickly explain the strategy of the proof, and we later show that the initial condition belongs to the desired space, i.e. $H^{11/4} (\ln H)^{-\beta}$. We also go through technical lemmas that we will need later to perform the proof of the blow-up.

\subsection{Explicit resolution and preliminary results}

We will first solve the equation using the characteristic method.

We consider the equation (\ref{modelequation}) and look at solutions of the form $u(t,x) = u_1(t,x_1)$. 

The equation in one space dimension can be factored as in the following.

\begin{equation}\label{dimone}
	\begin{split}
		\left( \left( \partial_t + \partial_{x_1} \right) + v \left( \partial_{x_1} - \partial_{t}\right) \right) \left( \partial_{x_1} -\partial_t \right) u =0 &\\
		u(0,x_1) =0, \hspace{0.4cm} \partial_t u(0,x_1) = -\chi(x_1)& 
	\end{split}
\end{equation}
where $v = (\partial_{x_1} - \partial_{t}) u$; which is equivalent when $v\neq 1$ to

\begin{equation}
	\left(\partial_t + \frac{1+v}{1-v} \partial_{x_1} \right) \left( \partial_{x_1} - \partial_t \right) u = 0.
\end{equation}

Now, this partial differential equation can be explicitly solved. Introducing $\phi$ such that

\begin{equation}\label{ceciestencoreungenrededef}
	\left\{
	\begin{aligned}
		&\phi(0,y) = y, \\
		&\partial_t \phi(t,y) = \frac{1+v(t,\phi(t,y))}{1-v(t,\phi(t,y)},
	\end{aligned}
	\right.
\end{equation}

we obtain

\begin{equation}
	\frac{\partial}{\partial t}\left(v(t,\phi(t,y)) \right) = 0\hspace{0.15cm} \Rightarrow \hspace{0.15cm} v(t,\phi(t,y)) = \chi(y) \hspace{0.2cm} \forall t.
\end{equation}

Now, this gives us an explicit formula for $\phi$,

\begin{equation}
	\phi(t,y) = y + t \frac{1+\chi(y)}{1-\chi(y)}.
\end{equation}

We will consider the Cauchy problem corresponding to the two following choices for the initial condition. First, we consider for $\alpha >0 \in \mathbb{R}$ (the conditions that $\alpha$ has to satisfy will appear later through the proof)

\begin{equation}\label{defchi4}
	v_0(x_1,x_2) = \chi(x_1) = - \int_{0}^{x_1} | \ln(s) |^\alpha ds,
\end{equation}

which will be the initial value that corresponds to an instantaneous blow-up of the solution. But we will work with the regularized initial condition 

\begin{equation}\label{defchi5}
	v_{0,\varepsilon} (x_1,x_2) = \chi_\varepsilon(x_1) = - \int_{0}^{x_1} \psi_\varepsilon(s) | \ln(s) |^\alpha ds,
\end{equation}

where 

\begin{equation}\label{defdefpsi}
	\left\{
	\begin{split}
		&1 > \psi_\varepsilon(x) > 0 \text{ for } \varepsilon/2 <x<\varepsilon, \\
		&\psi_\varepsilon(x) = 0 \text{ for } x< \varepsilon/2, \\
		&\psi_\varepsilon(x) = 1 \text{ for } x> \varepsilon \\
		&\exists C>0, |\psi_\varepsilon'(x)| \leq \frac{C}{\varepsilon}.
	\end{split}
	\right.
\end{equation}

This resolution holds whenever the solution depends only on $x_1$. We will consider initial conditions defined on $\mathbb{R}^2$ entirely, but that coincides with (\ref{defchi4}) or (\ref{defchi5}) on a set $\Omega$.

We choose 

\begin{equation}\label{condcondomega}
	\Omega_0 = \left\{ (x_1,x_2) | x_1 \geq 0, |x_2| \leq \frac{\sqrt{x_1}}{|\ln(x_1)|^\delta} \right\} \bigcap \left[0;\frac{1}{2}\right] \times \left[0;\frac{1}{2}\right],
\end{equation}

and we take $\Omega$ to be the largest domain of dependence (defined in Definition \ref{dependence}) for the metric whose inverse is given by
\begin{equation}\label{invmet}
	\sum_{i,j=0}^2 g^{ij} (t,x) \partial_{x_i} \partial_{x_j} = \partial_t^2 - \sum_{i=1}^2 \partial_{x_i}^2 - v(t,x_1)(\partial_t - \partial_{x_1})^2,
\end{equation}
and such that $\Omega\cap \{ t=0 \}= \Omega_0$. We correspondingly define $\Omega_t = \{ (x_1,x_2) | (t,x_1,x_2) \in \Omega \}$ and $a_t(x_1)$ to be the width of $\Omega_t$ at $x_1$.

\section{Introduction of the logarithmic perturbation of $H^s$ and related lemmas}

In this chapter, we will introduce logarithmic perturbations of the spaces $\dot H^{7/4}$ and $H^{7/4}$, namely $\dot H^{7/4} (\ln H)^{-\beta}$ and $H^{7/4} (\ln H)^{-\beta}$ that will contain functions slightly less (provided that $\beta\geq 0$) regular than $\dot H^{7/4}$ and $H^{7/4}$. We will show preliminary lemmas, and then we will show that we can find an extension of our function on $\mathbb{R}^2$ that has a uniformly bounded norm in this space. We define the space $\dot H^{s}(\ln H)^{-\beta}$ as the set of all functions such that the norm

\begin{multline}\label{logsobo}
	||f||_{\dot H^{s}(\ln H)^{-\beta}}^2 = ||\mathcal{F}(f)\cdot \frac{|\xi|^{s}}{\left( 1+|\ln(|\xi|)| \right)^\beta}||_{L^2}^2 \\
	= \int\int_{\xi_1,\xi_2\in \mathbb{R}^2} \left[ \frac{|\xi|^{s}}{\left( 1+|\ln(|\xi|)| \right)^\beta} \int\int_{x_1,x_2\in \mathbb{R}} e^{2i\pi \xi_1 x_1} e^{2 i \pi \xi_2 x_2} f(x_1,x_2) dx_1 dx_2 \right]^2 d\xi_1 d\xi_2,
\end{multline}

is finite. Similarly, we define the space $H^{7/4} (\ln H)^{-\beta}$ as the set of all functions such that the norm where $|\varepsilon|$ is replaced by $(1+|\varepsilon|^2)^{1/2}$ is finite.

We will first show that $\dot H^{s} \subseteq \dot H^{s}(\ln H)^{-\beta} \subseteq \dot H^{s-\lambda}$. This first lemma is for $L^1$ functions.

\begin{lemma}\label{dotHlog}
	Let $f$ be a function in $L^1(\mathbb{R}^2)$. For any nonnegative $s$ and $\beta$, for any small enough and positive $\lambda$, we have the following properties.
	\begin{itemize}
		\item[(i)] $||f||_{\dot H^{s}} < \infty \Rightarrow ||f||_{\dot H^{s}(\ln H)^{-\beta}} <\infty$,
		\item[(ii)] $||f||_{\dot  H^{s}(\ln H)^{-\beta}} <\infty  \Rightarrow ||f||_{\dot H^{s-\lambda}} < \infty$.
	\end{itemize}
\end{lemma}

\begin{proof}
	First, we note that because $f$ is in $L^1(\mathbb{R}^2)$, $\mathcal{F}(f)$ is globally bounded, indeed,
	
	\begin{multline}
		|\mathcal{F}(f)(\xi)| = \left| \int\int_{x_1,x_2\in \mathbb{R}^2} e^{2i\pi (x_1 \xi_1 + x_2 \xi_2)} f(x_1,x_2) dx_1 dx_2 \right| \leq \int\int_{x_1,x_2\in \mathbb{R}^2} |f(x_1,x_2)|  \\
		\leq ||f||_{L^1(\mathbb{R}^2)}.
	\end{multline}
	
	Now, we show $(i)$. Consider $f$ such that $||f||_{\dot H^{s}} < \infty$. Because $1+|\ln(|\xi|)| \geq 1$ for any $\xi$, we have 
	
	\begin{multline}
		||f||_{\dot H^{7/4}(\ln H)^{-\beta}}^2 = \int\int_{\xi_1,\xi_2\in \mathbb{R}^2} \frac{|\xi|^{2s}}{\left(1+|\ln(|\xi|)|^{\beta}\right)^2} \left(\mathcal{F}(f)(\xi_1,\xi_2)\right)^2 \\ 
		\leq  \int\int_{\xi_1,\xi_2\in \mathbb{R}^2} |\xi|^{2s} \left(\mathcal{F}(f)(\xi_1,\xi_2)\right)^2 =||f||_{\dot H^{s}}^2.
		\\
	\end{multline}
	
	Now, we show $(ii)$. Consider $f$ such that $||f||_{\dot H^{s}} < \infty$. Take $r(\lambda,\beta)$ such that $|\xi|>r$ implies $(1+|\ln(|\xi|)|^\beta)^2<|\xi|^{2\lambda}$. We then have
	
	\begin{multline}
		||f||_{\dot H^{s-\lambda}}^2 = \int\int_{\xi_1,\xi_2\in \mathbb{R}^2} \frac{|\xi|^{2s}}{|\xi|^{2\lambda}} \left(\mathcal{F}(f)(\xi_1,\xi_2)\right)^2 \\
		=  \int\int_{|\xi|<r(\lambda,\beta)} \frac{|\xi|^{2s}}{|\xi|^{2\lambda}} \left(\mathcal{F}(f)(\xi_1,\xi_2)\right)^2 + \int\int_{|\xi|>r(\lambda,\beta)} \frac{|\xi|^{2s}}{|\xi|^{2\lambda}} \left(\mathcal{F}(f)(\xi_1,\xi_2)\right)^2 \\
		\leq |B(0,r(\lambda,\beta))| \cdot |r(\lambda,\beta)|^{2s-2\lambda} ||\mathcal{F}(f)||_\infty^2 + \int\int_{|\xi|>r_0} \frac{|\xi|^{2s}}{(1+|\ln(|\xi|)|^\beta)^2} \left( \mathcal{F}(f)(\xi_1,\xi_2) \right)^2 \\
		\leq C(\lambda,\beta) \cdot ||\mathcal{F}(f)||_\infty + ||f||_{\dot H^{s}(\ln H)^{-\beta}} < \infty. 
	\end{multline}
	
\end{proof}

We now state a corollary, that is more relevant to our situation. 

\begin{lemma}\label{dotHlog2}
	Let $f \in L^2$.  
	Let $\beta$ and $s$ be two non-negative real numbers, $\lambda$ a small enough real number and $K$ a compact subset of $\mathbb{R}^2$. If $f$ is supported in $K$, then we have the two following:
	\begin{itemize}
		\item[(i)] $||f||_{\dot H^{s}} < \infty \Rightarrow ||f||_{\dot H^{s}(\ln H)^{-\beta}} <\infty$,
		\item[(ii)] $||f||_{\dot  H^{s}(\ln H)^{-\beta}} <\infty  \Rightarrow ||f||_{\dot H^{s-\lambda}} < \infty$.
	\end{itemize}
\end{lemma}

We now state one more lemma. Using this lemma, we will only have to compute the homogeneous logarithmically modified Sobolev norm as long as our functions are compactly supported. We will use this lemma, later on, to show that the initial condition we consider belongs to $H^{7/4} (\ln H)^{-\beta}$.

\begin{lemma}\label{celuila}
	Let $s$ be a non-negative real number, and $f\in L^1_{loc}$. If $f$ belongs to $ \dot H^{s} (\ln H)^{-\beta}$ and is compactly supported, then $f$ belongs to $H^s (\ln H)^{-\beta}$.
\end{lemma}

\begin{proof}
	Let $f\in \dot H^s (\ln H)^{-\beta}$ and supported in $K$, a compact subset of $\mathbb{R}^2$. Then, by lemma \ref{dotHlog2}, $f\in \dot H^{s/2}$. Hence, we have that $f \in L^2$. (The proof of this can for instance be found in \cite{bahouri2011fourier} p. 39).
	
	Now, we compute the non-homogeneous modified Sobolev norm.
	
	\begin{multline}
		\int\int_{\xi \in \mathbb{R}^2} \frac{(1+|\xi|^2)^s}{(1+|\ln(|\xi|)|^\beta)^2} |\hat f(\xi)|^2 \\
		= \int\int_{\xi \in B(0,1)} \frac{(1+|\xi|^2)^s}{(1+|\ln(|\xi|)|^\beta)^2} |\hat f(\xi)|^2 + \int\int_{\xi \in \mathbb{R}^2\setminus B(0,1)} \frac{(1+|\xi|^2)^s}{(1+|\ln(|\xi|)|^\beta)^2} |\hat f(\xi)|^2 \\
		\leq 2^s ||\hat f||_{L^2}^2 + \int\int_{\xi \in \mathbb{R}^2\setminus B(0,1)} \frac{(1+|\xi|^2)^s}{|\xi|^{2s}} \frac{|\xi|^{2s}}{(1+|\ln(|\xi|)|^\beta)^2} |\hat f(\xi)|^2 \\
		\leq 2^s \cdot C ||f||_{L^2} + 2^s \cdot C ||f||_{\dot H^s (\ln H)^{-\beta}}.
	\end{multline}
	
\end{proof}

Now, we will show a new lemma that we will use later. We express the Sobolev norm as a convolution-type integral. This type of integrals is widely used for the differentiation of fractional order. 

Used together with our previous lemma, it establishes a link between our logarithmically modified Sobolev spaces defined via Fourier transform and the fractional derivative.

\begin{lemma}\label{joachim} Let $\lambda$ be a small, nonnegative number. 
	Let $f: \mathbb{R}^2\to \mathbb{R}$ and $\omega \subseteq \mathbb{R}^2$ such that $f=0$ outside of $\omega$.
	
	Then, 
	
	\begin{equation}\label{joachimeq}
		\begin{split}
			||f||_{\dot H_{x_1}^{7/4-\lambda}(\mathbb{R}^2)}&= C  \int\int_{(x_1,x_2) \in \omega} \left( \frac{\partial^2 f}{\partial x_1^2} \right)(x_1,x_2)\\
			& \cdot  \int_{y|\hspace{0.1cm} (y,x_2) \in \omega} |x_1-y|^{-1/2+2\lambda} \left( \frac{\partial^2 f}{\partial x_1^2} \right)(y,x_2) dy dx_2 dx_1 \\
		\end{split}
	\end{equation}
	
\end{lemma}

\begin{proof}
	
	\begin{multline}
		||f||_{\dot H_{x_1}^{7/4-\lambda}(\mathbb{R}^2)}
		= ||\frac{\partial^2 f}{\partial x_1^2}||_{\dot H^{-1/4-\lambda}} \\
		= \int\int_{(x_1,x_2)\in\mathbb{R}^2} |\nabla_{x_1}^{-1/4-\lambda}| \left( \frac{\partial^2 f}{\partial x_1^2} \right)(x_1,x_2) \cdot |\nabla_{x_1}^{-1/4-\lambda}| \left( \frac{\partial^2 f}{\partial x_1^2} \right)(x_1,x_2) dx
		\\
		= \int\int_{(x_1,x_2)\in\mathbb{R}^2} |\nabla_{x_1}^{-1/2-2\lambda}| \left( \frac{\partial^2 f}{\partial x_1^2} \right)(x_1,x_2) \cdot \left( \frac{\partial^2 f}{\partial x_1^2} \right)(x_1,x_2) dx \\
		=_{(*)} C \int\int_{(x_1,x_2)\in\mathbb{R}^2} \left( \frac{\partial^2 f}{\partial x_1^2} \right)(x_1,x_2)\cdot \int_{y\in \mathbb{R}} \frac{\frac{\partial^2 f}{\partial x_1^2}(x_1,x_2)}{|x_1-y|^{1/2-2\lambda}} dx\\
		= C \int\int_{(x_1,x_2) \in \omega} \left( \frac{\partial^2 f}{\partial x_1^2} \right)(x_1,x_2)
		\cdot  \int_{y|\hspace{0.1cm} (y,x_2) \in \omega} |x_1-y|^{-1/2+2\lambda} \left( \frac{\partial^2 f}{\partial x_1^2} \right)(y,x_2) \\
	\end{multline}
	
	For $(*)$, we used that
	
	\begin{equation}\label{nabla1}
		\left( - \Delta\right)^{s/2} (f)(x) = \left( \left( 2\pi \left| \xi \right|\right)^s \widehat{f} (\xi) \right)^{\widecheck{\hspace{0.1cm}}} (x),
	\end{equation}
	
	and that for $n>s>0$,
	
	\begin{equation}\label{nabla2}
		\left( 2\pi \right)^{-s} \left( \left| \xi \right|^{-s} \right)^{\widecheck{\hspace{0.1cm}}}(x) = (2\pi)^{-s} \frac{\pi^{\frac{s}{2}}\Gamma(\frac{n-s}{2})}{\pi^{\frac{n-s}{2}}\Gamma(\frac{s}{2})} |x|^{s-n}.
	\end{equation}
	
	In our case, because we integrate only with respect to $x_1$, we obtain from (\ref{nabla1}) and (\ref{nabla2})
	
	\begin{equation}
		|\nabla_{x_1}|^{-1/2-2\lambda} f (x) = C(\lambda) \int_{y\in \mathbb{R}} \frac{f(y)}{|x_1-y|^{1/2-2\lambda}},
	\end{equation}
	
	where
	
	\begin{equation}
		C(\lambda)= (2\pi)^{-1/2-2\lambda} \pi^{2\lambda} \frac{\Gamma(\frac{1-\frac{1}{2}-2\lambda}{2})}{\Gamma(\frac{1/2+2\lambda}{2})}
	\end{equation}
	
\end{proof}

\section{Proof that $(\partial_t u)_{|t=0} \in H^{7/4}(\ln H)^{-\beta}$}\label{inisec}

We now introduce the main theorem of this chapter.

\begin{theorem}\label{inicondimportant}\label{thmextension}
	Let $\psi_\varepsilon$ be functions such that 
	
	\begin{equation}
		\left\{
		\begin{split} 
			&\psi_\varepsilon(x) = 0, \hspace{0.2cm} x\in [0,\varepsilon/2], \\
			&\psi_\varepsilon(x) = 1, \hspace{0.2cm} x\in [\varepsilon,\infty], \\
			&\forall x, ~ \psi_\varepsilon(x) \in [0,1], \\
			&|\partial_k \psi_\varepsilon(x)| \leq \frac{C(k)}{\varepsilon^k}.
		\end{split}
		\right.
	\end{equation}
	
	We consider $\alpha, \beta, \delta$ such that $2\alpha-2\beta-\delta < -1$. For a fixed $\beta>1/2$, it is possible to choose $\alpha>0$ and $\delta>0$ that satisfy this condition.
	
	With $\chi_\varepsilon:(x_1,x_2) \in \Omega \mapsto -\int_{0}^{x_1} \psi_\varepsilon(s) |\ln(y)|^\alpha dy$, there exists $h_\varepsilon: \mathbb{R}^2 \to \mathbb{R}$ such that $h_{\varepsilon,|\Omega_0}(x_1,x_2) = \chi_\varepsilon(x_1)$ and $||h_\varepsilon||_{H^{7/4}(\ln H)^{-\beta}} < \infty $. Moreover, the bound on the norm can be chosen to be independent of $\varepsilon$.
\end{theorem}

\begin{proof}
	First, we define $\chi_\varepsilon$ on $\mathbb{R}$ entirely by
	
	\begin{equation}
		\chi_\varepsilon:x_1 \in \mathbb{R} \mapsto \left\{
		\begin{split}
			&-\int_{0}^{x_1} \psi_\varepsilon(s) |\ln(y)|^\alpha dy, \text{ for } x_1>0, \\
			&0, \text{ for } x_1\leq 0.
		\end{split}
		\right.
	\end{equation}
	
	We consider a smooth function $\psi:\mathbb{R} \to [0,1]$ such that 
	
	\begin{equation}
		\left\{
		\begin{split} 
			&\psi(x) = 1, \hspace{0.2cm} x\in [0,1/4], \\
			&\psi(x) = 0, \hspace{0.2cm} x\geq 1/2. 
		\end{split}
		\right.
	\end{equation}
	
	$\psi$ is defined on $\mathbb{R}^-$ by setting $\psi(x) = \psi(-x)$ for $x<0$.

	Define 
	
	\begin{equation}
		h_\varepsilon(x_1,x_2) = \chi_\varepsilon(x_1) \cdot \psi\left(\frac{\left| \ln (x_1) \right|^\delta x_2}{\sqrt{x_1}} \right) \cdot \psi(x_1).
	\end{equation}
	
	We multiply $\chi$ by a cutoff function in $x_1$ and $x_2$ that respects geometry of $\Omega$, i.e. $\psi\left(\frac{\left|\ln (x_1) \right|^\delta x_2}{\sqrt{x_1}} \right) =0 $ when $x_2 \geq \frac{1}{2} \sqrt{x_1} |\ln (x_1)|^{-\delta}$ ; and we multiply $\chi$ by a simple cutoff function in $x_1$.
	
	Lastly, we consider a dyadic partition of unity, and $\lambda$ will denote dyadic numbers. Take a function $\zeta: \mathbb{R}\to [0,1]$ such that
	
	\begin{equation}
		\sum_{j\in \mathbb{Z}} \zeta_j(x_1) = \sum_{j\in\mathbb{Z}} \zeta\left( \frac{x_1}{2^j} \right) = \sum_\lambda \zeta\left(\frac{x_1}{\lambda}\right) =  1, ~ \forall x_1\in \mathbb{R},
	\end{equation}
	
	\begin{equation}
		Supp ~ \zeta_\lambda \subseteq [\frac{1}{4} \cdot 2^{j},4 \cdot 2^j] = [\frac{1}{4}\cdot \lambda,4 \cdot \lambda],
	\end{equation}
	
	\begin{equation}
		\forall x, ~\zeta_j(x) \in [0,1]. 
	\end{equation}
	
	Now, we define
	
	\begin{equation}
		h_{\lambda,\varepsilon}(x_1,x_2) = \zeta_\lambda (x_1) h_\varepsilon(x_1,x_2),
	\end{equation}
	
	and we have
	
	\begin{equation}
		h_\varepsilon = \sum_\lambda h_{\lambda,\varepsilon} = \sum_{\lambda\leq 2^{j_0}} h_{\lambda,\varepsilon}.
	\end{equation}
	
	In virtue of lemma \ref{celuila}, we only have to study the homogeneous modified Sobolev norm, as our function is compactly supported.
	
	Now, we will find an estimate for $||h_{\varepsilon,\lambda}||_{H^{7/4}(\ln H)^{-\beta}}$.
	
	First, we compute $\mathcal{F}\left( h_{\varepsilon,\lambda} (\cdot,\cdot) \right)(\xi_1,\xi_2)$.
	
	\begin{multline}
		\mathcal{F}\left( h_{\varepsilon,\lambda} (\cdot,\cdot) \right)(\xi_1,\xi_2) = \int_{x_1} \int_{x_2} e^{-2i \pi x_1 \xi_1} e^{-2i \pi x_2 \xi_2}  \psi(x_1)\zeta_{\lambda}(x_1) \chi_\varepsilon(x_1)  \psi\left( \frac{|\ln(x_1)|^\delta}{\sqrt{x_1}} x_2 \right) \\
		=\int_{x_1} e^{-2i\pi xi_1}  \psi(x_1)\zeta_{\lambda}(x_1) \chi_\varepsilon(x_1)  \int_{x_2} e^{-2i\pi \xi_2} \psi\left( \frac{|\ln(x_1)|^\delta}{\sqrt{x_1}} x_2 \right) \\
		= \int_{x_1} e^{-2i \pi xi_1 } \psi(x_1) \zeta_\lambda(x_1) \chi_\varepsilon(x_1) \left[ \frac{\sqrt{x_1}}{|\ln(x_1)|^\delta} \cdot \mathcal{F} \left(\psi\right) \left( \frac{\sqrt{x_1}}{|\ln(x_1)|^\delta} \xi_2\right) \right] \\
		= \int_{x_1=\frac{\lambda}{4}}^{4\lambda} e^{-2i \pi xi_1 } \psi(x_1) \zeta_\lambda(x_1) \chi_\varepsilon(x_1) \left[ \frac{\sqrt{x_1}}{|\ln(x_1)|^\delta} \cdot \mathcal{F} \left(\psi\right) \left( \frac{\sqrt{x_1}}{|\ln(x_1)|^\delta} \xi_2\right) \right].
	\end{multline}
	
	Because $\frac{1-C}{C} \rightarrow_{C\rightarrow 0^+} \infty$, we can chose $C_0>0$ such that $\frac{1-C_0}{C_0}  \geq \frac{\alpha}{|\ln \frac{1}{2}|^\alpha}$.
	
	Now, we have for $ \varepsilon/2 \leq y \leq 1/2$, and $\alpha \leq 1$
	
	\begin{equation}
		|\ln (y)|^\alpha - \alpha |\ln (y)|^{\alpha-1} \geq C_0 |\ln(y)|^\alpha.
	\end{equation}
	
	Hence, we have for $x_1\in [\lambda/4,4 \lambda]$, 
	
	\begin{multline}
		\left| \chi_\varepsilon(x_1) \right| = \left| \int_{s=\varepsilon/2}^{x_1} \psi_\varepsilon(s) |\ln(s)|^\alpha \right| \\
		\leq  \frac{1}{C_0} \left| \int_{s=\varepsilon/2}^{x_1} |\ln(s)|^\alpha - \alpha |\ln(s)|^{\alpha-1}  \right| \leq  C \lambda |\ln(\lambda)|^\alpha.
	\end{multline}
	
	Now, we look at
	
	\begin{equation}
		||h_{\lambda,\varepsilon}||_{\dot H^{7/4}(\ln H)^{-\beta}}^2 = \int_{\xi_1,\xi_2} \left( \frac{|\xi|^{7/4}}{(1+|\ln(|\xi|)|)^\beta} \cdot \mathcal{F}(h_{\lambda,\varepsilon})(\xi_1,\xi_2) \right)^2.
	\end{equation}
	
	We will use the fact that $\mathcal{F}(\phi)\left( \frac{\sqrt{x_1} \xi_2}{|\ln(x_1)|^\delta} \right)$ is rapidly decreasing when $\xi_2 >> \frac{|\ln(\lambda)|^\delta}{\sqrt{\lambda}}$ and that $\mathcal{F}(\zeta\lambda)(\xi_1) = \mathcal{F}(\zeta(\frac{\cdot}{\lambda}))(\xi_1)$ is rapidly decreasing when $\xi_1>>\frac{1}{\lambda}$ to essentially reduce the integration domain to $[0,\frac{1}{\lambda}]\times [0,\frac{|\ln(\lambda)|^\delta}{\sqrt{\lambda}}]$.
	
	First, we compute the following integral
	
	\begin{equation}
		\int_{\xi_1\leq \lambda^{-1}}\int_{\xi_2 \leq \frac{|\ln(\lambda)|^\delta}{\sqrt{\lambda}}} \left( \frac{|\xi|^{7/4}}{(1+|\ln(|\xi|)|)^\beta} \cdot \mathcal{F}(h_{\lambda,\varepsilon})(\xi_1,\xi_2) \right)^2.
	\end{equation}
	
	We have that
	
	\begin{equation}
		\left| \mathcal{F}(h_{\lambda,\varepsilon})(\xi_1,\xi_2) \right| \leq C \left| \int_{\lambda/4}^{4\lambda} \lambda |\ln(\lambda)|^\alpha \frac{\sqrt{x_1}}{|\ln(x_1)|^\delta} \right| \leq C \lambda^{5/2} |\ln(\lambda)|^{\alpha-\delta}.
	\end{equation}
	
	And so, we get
	
	\begin{multline}
		\left| \int_{\xi_1\leq \lambda^{-1}}\int_{\xi_2 \leq \frac{|\ln(\lambda)|^\delta}{\sqrt{\lambda}}} \left( \frac{|\xi|^{7/4}}{(1+|\ln(|\xi|)|)^\beta} \cdot \mathcal{F}(h_{\lambda,\varepsilon})(\xi_1,\xi_2) \right)^2 \right| \\
		\leq C \lambda^{-1} \frac{|\ln(\lambda)|^\delta}{\sqrt{\lambda}} \frac{\lambda^{-7/2}}{|\ln(\lambda)|^{2\beta}} \lambda^5 |\ln(\lambda)|^{2\alpha - 2\delta} = C |\ln(\lambda)|^{2\alpha-2\beta - \delta} .
	\end{multline}
	
	Now, we will make a precise argument to justify that integrating over the whole space $\mathbb{R}^2$ does not give a bigger term in $\lambda$.
	
	First, we look at 
	
	\begin{equation}
		\int_{\xi_2 \leq \frac{|\ln(\lambda)|^\delta}{\sqrt{\lambda}}}  \int_{\xi_1\geq \lambda^{-1}} \left( \frac{|\xi|^{7/4}}{(1+|\ln(|\xi|)|)^\beta} \cdot \mathcal{F}(h_{\lambda,\varepsilon})(\xi_1,\xi_2) \right)^2.
	\end{equation}
	
	We can write
	
	\begin{multline}
		\int_{\lambda/4}^{4\lambda} e^{2i \pi \xi_1 x_1} \chi_\varepsilon(x_1) \zeta_\lambda(x_1) \frac{\sqrt{x_1}}{|\ln(x_1)|^\delta} \mathcal{F}(\phi)\left( \frac{\sqrt{x_1}}{|\ln(x_1)|^\delta} \right)\\
		= C \frac{1}{\xi_1} \int_{\lambda/4}^{4\lambda} e^{2i \pi \xi_1 x_1} \Big[ \chi_\varepsilon'(x_1) \zeta_\lambda(x_1) \frac{\sqrt{x_1}}{|\ln(x_1)|^\delta} \mathcal{F}(\phi)\left( \frac{\sqrt{x_1}}{|\ln(x_1)|^\delta} \right) \\
		+ \chi_\varepsilon(x_1) \zeta_\lambda'(x_1) \frac{\sqrt{x_1}}{|\ln(x_1)|^\delta} \mathcal{F}(\phi)\left( \frac{\sqrt{x_1}}{|\ln(x_1)|^\delta} \right) \\
		+\chi_\varepsilon(x_1) \zeta_\lambda(x_1) \frac{1}{2\sqrt{x_1} |\ln(x_1)|^\delta} \mathcal{F}(\phi)\left( \frac{\sqrt{x_1}}{|\ln(x_1)|^\delta} \right) \\
		+\chi_\varepsilon(x_1) \zeta_\lambda(x_1) \frac{+\delta}{x} \frac{\sqrt{x_1}}{ |\ln(x_1)|^{\delta+1}} \mathcal{F}(\phi)\left( \frac{\sqrt{x_1}}{|\ln(x_1)|^\delta} \right) \\
		+\chi_\varepsilon(x_1) \zeta_\lambda(x_1) \left( \frac{1}{2 |\ln(x_1)|^{2\delta}} + \frac{\delta}{|\ln(x_1)|^{2\delta+1}} \right) \left( \mathcal{F}(\phi) \right)' \left( \frac{\sqrt{x_1}}{|\ln(x_1)|^\delta} \right) \Big].
	\end{multline}
	
	Now, the $\frac{1}{\xi_1}$ we gain is going to be smaller than $\lambda$ on our considered set. So now, we show that we lose at most $\lambda$ when differentiating the involved functions. When we will integrate $\frac{1}{\xi_1^2}$ over the set $\xi_1\geq 1/\lambda$, we will multiply by $\frac{1}{\lambda}$ which is not worse than the $\frac{1}{\lambda}$ we had in the first estimate because of the size of the set. 
	
	First,
	
	\begin{equation}\label{encoreun1}
		\chi_\varepsilon'(x_1) \leq C |\ln(\lambda)|^\alpha,
	\end{equation}
	
	\begin{equation}\label{encoreun2}
		\zeta_\lambda'(x_1) = \frac{\partial}{\partial_{x_1}} \left( \zeta \left( \frac{x_1}{\lambda} \right) \right) (x_1) = \frac{1}{\lambda} \zeta'\left(\frac{x_1}{\lambda} \right) \leq \frac{C}{\lambda}.
	\end{equation}
	
	So, we have
	
	\begin{equation}
		\left| \chi_\varepsilon'(x_1) \zeta_\lambda(x_1) \frac{\sqrt{x_1}}{|\ln(x_1)|^\delta} \mathcal{F}(\phi)\left( \frac{\sqrt{x_1}}{|\ln(x_1)|^\delta} \right) \right| \leq C |\ln(\lambda)|^{\alpha-\delta} \sqrt{\lambda},
	\end{equation}

	\begin{equation}
		\left| \chi_\varepsilon(x_1) \zeta_\lambda'(x_1) \frac{\sqrt{x_1}}{|\ln(x_1)|^\delta} \mathcal{F}(\phi)\left( \frac{\sqrt{x_1}}{|\ln(x_1)|^\delta} \right) \right| \leq C |\ln(\lambda)|^{\alpha-\delta} \sqrt{\lambda},
	\end{equation}
	
	\begin{equation}
		\left| \chi_\varepsilon(x_1) \zeta_\lambda(x_1) \frac{1}{2\sqrt{x_1} |\ln(x_1)|^\delta} \mathcal{F}(\phi)\left( \frac{\sqrt{x_1}}{|\ln(x_1)|^\delta} \right) \right| \leq C |\ln(\lambda)|^{\alpha-\delta} \sqrt{\lambda},
	\end{equation}
	
	\begin{equation}
		\left| \chi_\varepsilon(x_1) \zeta_\lambda(x_1) \frac{\delta}{\sqrt{x_1} |\ln(x_1)|^{\delta+1}} \mathcal{F}(\phi)\left( \frac{\sqrt{x_1}}{|\ln(x_1)|^\delta} \right) \right| \leq C |\ln(\lambda)|^{\alpha-\delta} \sqrt{\lambda},
	\end{equation}
	
	\begin{multline}
		\Big| \chi_\varepsilon(x_1) \zeta_\lambda(x_1) \left( \frac{1}{2 |\ln(x_1)|^{2\delta}} + \frac{\delta}{|\ln(x_1)|^{2\delta+1}} \right) \\
		\cdot \left( \mathcal{F}(\phi) \right)' \left( \frac{\sqrt{x_1}}{|\ln(x_1)|^\delta} \right) \Big| \leq C |\ln(\lambda)|^{\alpha-\delta} \sqrt{\lambda}.
	\end{multline}
	
	By induction and Leibniz differentiation formula, we quickly obtain that
	
	\begin{equation}
		\frac{\partial^k}{\partial_{x_1}^k} \left( \chi_\varepsilon(x_1) \zeta_\lambda(x_1) \frac{\sqrt{x_1}}{|\ln(x_1)|^\delta} \mathcal{F}(\phi)\left( \frac{\sqrt{x_1}}{|\ln(x_1)|^\delta} \right) \right) \leq C_k |\ln(\lambda)|^{\alpha-\delta} \lambda^{\frac{3}{2} - k}.
	\end{equation}
	
	Hence, we have
	
	\begin{multline}\label{xi1}
		\left( \int_{\lambda/4}^{4\lambda} e^{2i \pi \xi_1 x_1} \chi_\varepsilon(x_1) \zeta_\lambda(x_1) \frac{\sqrt{x_1}}{|\ln(x_1)|^\delta} \mathcal{F}(\phi)\left( \frac{\sqrt{x_1}}{|\ln(x_1)|^\delta} \right) \right| \\
		\leq \frac{C_5}{\xi_1^5} |\ln(\lambda)|^{\alpha-\delta} \lambda^{-7/2} \cdot \lambda
	\end{multline}
	
	Now, we obtain that
	
	\begin{multline}
		\int_{\xi_1\geq \lambda^{-1}} \left( \frac{|\xi|^{7/4}}{(1+|\ln(|\xi|)|)^\beta} \cdot \mathcal{F}(h_{\lambda,\varepsilon})(\xi_1,\xi_2) \right)^2 \\
		\leq \int_{\xi_1\geq \lambda^{-1}} \frac{C_5}{\xi_1^{10} } |\ln(\lambda)|^{2\alpha-2\delta} \lambda^{-7} \frac{|\xi|^{7/2}}{|\ln(|\xi|)|^{2\beta}} \lambda^2  \\
		\leq \frac{C}{\lambda^{5}} \lambda^{11/2} |\ln(\lambda)|^{2\alpha-2\delta-2\beta} = C \lambda^{1/2} |\ln(\lambda)|^{2\alpha-2\beta-2\delta}. 
	\end{multline}
	
	Lastly, we hence obtain
	
	\begin{multline}
		\int_{\xi_2=0}^{\frac{|\ln(\lambda)|^\delta}{\sqrt{\lambda}}} \int_{\xi_1} \int_{\lambda/4}^{4\lambda} e^{2i \pi \xi_1 x_1} \chi_\varepsilon(x_1) \zeta_\lambda(x_1) \frac{\sqrt{x_1}}{|\ln(x_1)|^\delta} \mathcal{F}(\phi)\left( \frac{\sqrt{x_1}}{|\ln(x_1)|^\delta} \right) \leq C |\ln(\lambda)|^{2\alpha-2\beta-\delta}.
	\end{multline}
	
	The computations for $\xi_2 \geq  2 \frac{|\ln(\lambda)|^\delta}{\sqrt{\lambda}}$ are made using similar techniques, and lead to the same result.
	
	Now, we have that 
	
	\begin{equation}
		||\frac{|\xi|^{7/4}}{(1+|\ln(|\xi|)|)^\beta}\mathcal{F}(\zeta_\lambda \cdot v_0)(\xi_1,\xi_2)||_{L^2(\xi_1,\xi_2)}^2 \leq C |\ln(\lambda)|^{2\alpha-2\beta-\delta}. 
	\end{equation}
	
	Taking the sum over the dyadic numbers $\lambda=2^{-k}$, with $2\alpha-2\beta-\delta<-1$, we obtain that $g_{\varepsilon} \in H^{7/4}(\ln H)^{-\beta}$. Also, because the constant $C$ does not depend on $\varepsilon$, we obtain that the Sobolev norms of $g_{\varepsilon}$, for $\varepsilon\in (0,-1]$, are uniformly bounded.

\end{proof}

\section{Lower bound on the width of the domain near the singularity}

Lastly, we will use the following result that gives an estimation of $a_t(\phi(t,y))$ for small values of $y$. Here we only need the fact that the width obtained at the singularity is strictly bigger than zero. 

Using this result, we will be able to use a cutoff near the singularity in the next chapter. For this lemma, we consider the initial condition defined by (\ref{defchi5}) and with a cutoff as defined in theorem \ref{inicondimportant}. We consider $t_\varepsilon$ to be the first value such that there exists $\nu_\varepsilon$ such that $\phi_y(t_\varepsilon,\nu_\varepsilon)=0$. The fact that $t_\varepsilon$ and $\nu_\varepsilon$ exist and $t_\varepsilon= O(\frac{1}{|\ln(\varepsilon)|^\alpha})$ will be shown in chapter \ref{lesepsilons}.

\begin{proposition}\label{propdomainepro}
	Let $\nu_\varepsilon \in \mathbb{R}\setminus \{0\} $, and consider an interval of $t$ of the form $J_\varepsilon = [0,C \frac{1}{|\ln(\varepsilon)|^\delta}[$. Then there exists a positive constant $r$ such that
	
	\begin{equation}\label{1}
		a_0(\phi(0,y)) \sim_{y=0} 2 \sqrt{2} \frac{ \sqrt{y}}{|\ln(y)|^\delta},
	\end{equation}
	
	and for $t\in J_\varepsilon$,
	
	\begin{equation}\label{35lol}
		B\left((t,\phi(t,\nu_\varepsilon),0),r \right) = \{ (t,x_1,x_2) \in \mathbb{R}^3 | \sqrt{ (\phi(t,\nu_\varepsilon)-x_1)^2 + x_2^2 } \leq r \} \subset \Omega_t,
	\end{equation}
	
	provided that the condition 
	
	\begin{equation}\label{condalphadelta}
		\alpha > 2 \delta
	\end{equation}
	
	is satisfied.
	
\end{proposition}

\begin{proof}(of (\ref{1}))
	(trivial) Because $\phi(0,y) = y$, 
	
	\begin{equation}
		a_0(\phi(0,y)) = a_0(y) = \int_{|x_2| \leq \sqrt{2}\frac{ \sqrt{y}}{|\ln(y)|^\delta}} dx_2 = 2 \sqrt{2} \frac{ \sqrt{y}}{|\ln(y)|^\delta}
	\end{equation}
	
	Let us now prove \ref{35lol}. We assume that $|v| < 1/100$, and that $t\leq C \frac{1}{|\ln(\varepsilon)|^\alpha}$ (This will be achieved whenever $\varepsilon$ is small enough, which means we will only consider times $t \leq  C \frac{1}{|\ln(\varepsilon)|^\alpha}$ very small, see part 3 for further details.)
	
	We will distinguish three cases, first, we consider curves whose starting point has an abscissa strictly bigger than $x_0 = \frac{\varepsilon}{4}$.
	
	It follows from definition \ref{dependence} that $(t',x') \in \Omega$ if and only if $(t',x_1') \in \Omega^1$ and all Lipschitz continuous curves from $(t',x')$ that satisfy (\ref{condlight}) intersect the hyperplane $t=0$ in the set $\{ x | |x_2| \leq \frac{\sqrt{2x_1}}{|\ln(x_1)|^\delta} \}$. 
	
	Now, let $(t(s),x_1(s),x_2(s))$ be a Lipschitz continuous curve parameterized so that $t(s) + x_1(s) = s$. Note $q(s) = x_1(s) - t(s)$. Note that (\ref{condlight}) is equivalent to (using the fact that $\frac{dt(s)}{ds}+\frac{dx_1(s)}{ds}=1$),

	\begin{equation}
		R(s) \leq v(t(s),x_1(s)) - \frac{d q(s)}{ds},
	\end{equation}
	where $R(s) = \left(\frac{dx_2(s)}{ds}\right)^2$.
	
	Now, using this set of new variables $s= x_1 +t$, $q = x_1-t$ and $U(s,q) = u((s-q)/2,(s+q)/2)$, (\ref{dimone}) becomes
	
	\begin{equation}
		\left\{
		\begin{split}
			(\partial_s + V(s,q) \partial_q) \partial_q U(s,q) &= 0, \hspace{0.2cm} V(s,q) = 2 \partial_q U(s,q),\\
			U(y,y) &= 0, \hspace{0.2cm} U_q(y,y) = \frac{1}{2}\chi(y)
		\end{split}
		\right.
	\end{equation}
	
	The characteristics are given by $s = constant$ and $q =h(s,y)$ with 
	\begin{equation}
		\frac{d}{ds}h(s,y) = V(s,h(s,y)), \hspace{0.2cm} h(y,y)= y.
	\end{equation}
	Thus, $s \mapsto V(s,h(s,y))$ is constant on the curve and is equal to $\chi(y)$. These are the key ingredients to make this proof. The rest of the proof for this case is similar to the one that is done in \cite{lindblad1998counterexamples}. The reason why we need the hypothesis on $\alpha$ and $\delta$, and also we are not able to conclude with this method in the case $x_0 \leq \frac{\varepsilon}{4}$ comes from the fact that the domain is slightly less "wide" near $x=0$. For the detailed proof, see appendix \cite{appendixA}.

	Now, we study the case where the starting abscissa is smaller than $x_0=\varepsilon/4$.

	We will consider the curve $\mathcal{C}= \{ (x_1,x_2) | x_2 = \pm \frac{\sqrt{2x_1} }{|\ln(x_2)|^\delta}, x_0 \leq x_1<  t/2 \}$.
	
	The distance between the curve $\mathcal{C}$ and the point $(\phi(t_\varepsilon,\nu_\varepsilon),0)= (t,0)$ is given by
	\begin{equation}
		f(y) = \sqrt{ (\phi(t_\varepsilon,\nu_\varepsilon)-y)^2 + \frac{2y}{|\ln(y)|^{2 \delta} } }.
	\end{equation}
	
	Now, because for any $t<t_\varepsilon$, for any $y$, $\partial_y \phi(t,y) \neq 0$, and $\partial_y \phi(t,0)>0$, we obtain that for any $t<t_\varepsilon$, for any $y$, $\partial_y \phi(t,y) \geq 0$. Hence, we obtain the following inequality (because $\nu_\varepsilon \geq \varepsilon/2$)
	
	\begin{equation}
		\phi(t,\nu_\varepsilon ) \geq \phi(t,\varepsilon/2) = \frac{\varepsilon}{2} + t \frac{1+\chi_\varepsilon(\varepsilon/2)}{1-\chi_\varepsilon(\varepsilon/2)} = \frac{\varepsilon}{2} + t,
	\end{equation}
	
	and so we get $\phi(t_\varepsilon,\nu_\varepsilon) \geq \frac{\varepsilon}{2} + t_\varepsilon$. Now, because $y \leq \varepsilon/4$, we obtain that
	
	\begin{equation}
		f(y) \geq |\phi(t_\varepsilon,\nu_\varepsilon) - y| \geq \frac{\varepsilon}{4} + t.
	\end{equation}
	
	Let us rewrite the condition (\ref{dependence}). The metric $(g^{i,j})$ is given by (\ref{metricg}), which means that the inverse is given by
	\begin{equation}
		\left( g_{i,j}\right) =
		\begin{bmatrix}
			1+v & v & 0  \\
			v & -1+v & 0  \\
			0 & 0 & -1  \\
		\end{bmatrix}.
	\end{equation}
	
	Now, (\ref{dependence}) becomes
	\begin{equation}
		(1+v) \cdot 1 + 2 v \left( \frac{\partial x_1}{\partial_t} \right) + (-1+v) \left( \frac{\partial x_1}{\partial t} \right)^2 - \left( \frac{\partial  x_2}{\partial t} \right)^2 \geq 0,
	\end{equation}
	
	or after a few steps,
	
	\begin{equation}\label{ellipse}
		(1-v)^2 \left( \frac{\partial x_1}{\partial t}-\frac{v}{1-v} \right)^2 + (1-v) \left( \frac{\partial x_2}{\partial t} \right)^2 \leq 1.
	\end{equation}

	Call $E$ the ellipse given by (\ref{ellipse}), and $C$ the circle of center $(0,0)$ and radius $1$. We will show that we are in the situation depicted in figure \ref{piclol}.
	
	\begin{figure}[h]
		\centering
			\caption{}\label{piclol}
			\begin{tikzpicture}[scale=1.5]
				\draw[help lines, color=gray!30, dashed] (-1.4,-1.4) grid (1.4,1.4);
				\node[above] at (0.6,0) {\small $\frac{1}{1-v}$};
				\node[left] at (0.05,0.7) {\small $\frac{1}{\sqrt{1-v}}$};
				\node[left,below] at (-0.4,-0.1) {\small $(\frac{v}{1-v},0)$};
				\draw (-0.25,-0.05) -- (-0.25,0.05);
				\draw (-0.11,0.89) -- (-0.11,0.99);
				\node at (0.45,-0.55) {$E$};
				\node at (0.75,-0.55) {$C$};
				\draw[->,ultra thick] (-1.25,0)--(1.25,0) node[right]{\small $\frac{\partial x_1}{\partial t}$};
				\draw[->,ultra thick] (0,-1.25)--(0,1.25) node[above]{\small $\frac{\partial x_2}{\partial t}$};
				\draw (-0.11,0) ellipse (0.88cm and 0.94cm);
				\draw (0,0) circle (1cm);
			\end{tikzpicture}
	\end{figure}
	
	Let us prove that $E$ is in fact included in $C$.
	Let us compute $E\cap C$,

	\begin{equation}\label{inter}
		\left\{
		\begin{split}
			(1-v)^2 \left(x-\frac{v}{1-v}\right)^2 + (1-v) \left( y \right)^2 = 1 \\
			x^2 + y^2 = 1
		\end{split}
		\right.
	\end{equation}
	
	Assuming by symmetry $y\geq 0$,  we obtain $(1-v)^2 (x - \frac{v}{1-v})^2 + (1-v) (1-x^2) = 1 $. Now, the discriminant of this equation in $x$ is 
	\begin{equation}
		\Delta = \frac{4 v^2}{(1-v)^2} - 4 \frac{v}{1-v} \frac{v}{1-v} = 0,
	\end{equation}
	and the only solution we find is $(\frac{\frac{2v}{1-v}}{\frac{-2v}{1-v}},0) = (-1,0)$. This means that the ellipse $E$ is included in the circle $C$. This means that a curve satisfying (\ref{dependence}), satisfies 
	\begin{equation}
		\left( \frac{\partial x_1}{\partial t}\right)^2 + \left( \frac{\partial x_2}{\partial t}\right)^2 \leq 1
	\end{equation}
	and hence, any curve satisfying (\ref{dependence}) also satisfies 
	\begin{equation}\label{ceci}
		||(x_1(t),x_2(t)) - (x_1(0),x_2(0))||_2 \leq t.
	\end{equation}
	
	Because $d(\left(x_1(t),x_2(t)\right),C) \geq t+\frac{\varepsilon}{4}$, there exists a positive number $\delta_2$ such that a ball of radius $\delta_2$ and centered in $(\phi(t,\nu_\varepsilon),0)$ cannot be reached.
	
\end{proof}

\textbf{Reducing the domain to $x_1 \leq \frac{1}{|\ln(\varepsilon)|^{\alpha/2}}$.}

In this chapter, we additionally multiply our initial condition by a cutoff in $x_1$. The goal is to cut the function way after the point around which the phenomenon occurs, but to still have a domain that becomes as small as we want when $\varepsilon \rightarrow 0$. We will need the function to still be in $H^{7/4}(\ln H)^{-\beta}$, and we also need that the point $(\nu_\varepsilon,0)$ is in the interior of the domain of dependence at the time $t=t_\varepsilon$. 

For the second part of the requirements, a rough estimate is to notice that the speed of the information in our problem is at most $1$, and $t_\varepsilon$ satisfies $t_\varepsilon \leq \frac{1}{|\ln(\varepsilon)|^\alpha}$. This means that if the cutoff modifies the function only for $x_1 \geq \frac{1}{|\ln (\varepsilon)|^{\alpha/2}}$, thanks to (\ref{ceci}), we obtain our desired result for $\varepsilon$ small enough.

Hence, we consider a cutoff of the form (with the previous notations) $l_\varepsilon(x_1) = l\left(x_1 \cdot |\ln(\varepsilon)|^{\alpha/2}\right)$. We do not give more details here because there is no difficulty on this side in $x$.

\section{Blow-up of the solution to the regularized problem in $H^{11/4-\lambda}(\mathbb{R}^2)$.}\label{exavant}

In this chapter, we consider a regularized version of the initial condition, so the equation is well-posed and the method we use to compute its expression is sound. We then construct a counterexample using the statements we have made.

\section{Strategy}

In this chapter, we will consider the initial condition $g_\varepsilon$ whose existence and definition are provided in theorem \ref{inicondimportant} (it was denoted by $h_\varepsilon$ in the proof of the theorem).

\begin{equation}
	g_\varepsilon(x_1,x_2) = \chi_\varepsilon(x_1) \cdot \psi \left(\frac{ \left| \ln (x_1)\right|^\delta x_2}{\sqrt{x_1}} \right) \cdot \psi(x_1),
\end{equation}

where

\begin{equation}
	\begin{split}
		&\chi_\varepsilon(x) =-\int_{0}^{x} \psi_\varepsilon(s) \left| \ln(s) \right|^\alpha ds,\\
	\end{split}
\end{equation}

First, we recall that the initial condition is in $\dot H^{7/4}(\ln H)^{-\beta}$, with a norm that can be bounded uniformly with respect to $\varepsilon$. Also, we define 

\begin{equation}
	\kappa(x_1,x_2) = \psi \left(\frac{ \left| \ln (x_1)\right|^\delta x_2}{\sqrt{x_1}} \right) \cdot \psi(x_1).
\end{equation}

The Cauchy problem
\begin{equation}\label{modelequationepsilon}
	\left\{
	\begin{split}
		\Box u &= (D u) D^{2} u,  \\
		(u,\partial_t u)_{|t=0} &= (0,-g_\varepsilon),
	\end{split}
	\right.
\end{equation}

is now well-posed on some interval of the form $[0,t_\varepsilon^1[$. Then, we define some time $t_\varepsilon$ for which we start to observe the concentration of the characteristics described in the first chapter, for some point $\nu_\varepsilon$. By this, we mean that $\phi_{\varepsilon,y} (t_\varepsilon,\nu_\varepsilon) = 0$. The Sobolev norm of the solution will be proven to be unbounded as $t\rightarrow t_\varepsilon$ ; besides the time $t_\varepsilon$ is going to $0$ as $\varepsilon \rightarrow 0$. 
Next, using a scaling argument, we put together a sequence of these solutions for which the lifespan is going to $0$, and such that total initial Sobolev norm is still finite. 
For any time $t>0$, $t$ will be beyond the lifespan of one of those solutions, thus leading to an infinite Sobolev norm.

\section{Blow up when $t \rightarrow t_\varepsilon$ and control of $t_\varepsilon$ with respect to $\varepsilon$.}\label{lesepsilons}

Let $u_\varepsilon$ be the solution for $t < t_\varepsilon$ of the Cauchy problem:

\begin{equation}\label{Cauchyjo}
	\left\{
	\begin{split}
		\Box u_\varepsilon &= v_\varepsilon D v_\varepsilon,  \\
		(u,\partial_t u)_{|t=0} &= (0,-g_\varepsilon),
	\end{split}
	\right.
\end{equation}

with $v_\varepsilon = D u_\varepsilon$. We will prove the following theorem:

\begin{theorem}\label{conj}
	Let $u_\varepsilon$ be a solution of (\ref{Cauchyjo}), and $\delta_\varepsilon>0$ (conditions on $\delta_\varepsilon$ will be specified later).
	
	Let $\psi_{\varepsilon}^1: \mathbb{R} \to \mathbb{R}$ be a $C^\infty$ function satisfying
	
	\begin{equation}
		\left\{
		\begin{split}
			&\psi_\varepsilon^1(x) = 1 \text{ for } \phi(t_\varepsilon,\nu_\varepsilon) - \delta_\varepsilon <x<\phi(t_\varepsilon,\nu_\varepsilon) + \delta_\varepsilon \\
			&\psi_\varepsilon^1(x) = 0 \text{ for } \phi(t_\varepsilon,\nu_\varepsilon) + 2 \delta_\varepsilon <x \text{ or } x<\phi(t_\varepsilon,\nu_\varepsilon) - 2 \delta_\varepsilon \\
			&0 < \psi_\varepsilon^1(x) < 1 \text{ elsewhere},
		\end{split}
		\right.
	\end{equation}
	
	and $\psi_{\varepsilon}^2: \mathbb{R} \to \mathbb{R}$ be a $C^\infty$ function satisfying
	
	\begin{equation}
		\left\{
		\begin{split}
			&\psi_\varepsilon^2(x) = 1 \text{ for } - \delta_\varepsilon <x< \delta_\varepsilon \\
			&\psi_\varepsilon^2(x) = 0 \text{ for }  2 \delta_\varepsilon <x \text{ or } x< - 2 \delta_\varepsilon \\
			&0 < \psi_\varepsilon^2(x) < 1 \text{ elsewhere},
		\end{split}
		\right.
	\end{equation}
	
	so that $h_\varepsilon:(t,x_1,x_2) \mapsto v_\varepsilon(t,x_1,x_2) \psi_\varepsilon^1(x_1) \psi_\varepsilon^2(x_2)$ is localized in a square of width $4 \delta_\varepsilon$, cut in half by $x_1 = \phi_\varepsilon(t_\varepsilon,\nu_\varepsilon)$; and such that $h_\varepsilon=v_\varepsilon$ in a square of width $2 \delta_\varepsilon$, cut in half by $x_1 = \nu_{t,\varepsilon}$.
	
	We have, for any $\lambda >0$ small enough,
	
	\begin{equation}
		||h_\varepsilon(t)||_{H^{7/4-\lambda}_{x_1}} \rightarrow \infty \hspace{0.2cm} \text{ as } t\rightarrow t_\varepsilon.
	\end{equation}
	
\end{theorem}

\emph{Preliminary work}

Now, the computations made in the first part are still valid in $\Omega$, the domain of dependence such that $\Omega \cap \{t = 0 \} = \{ (x_1,x_2) | \chi_\varepsilon \kappa = \chi \}$. We also only consider the domain in time before the blow-up. Or, with $\phi_\varepsilon$ computed as previously,

\begin{equation}
	\Omega = \{ (x_1,x_2,t) | \phi_\varepsilon(t,\varepsilon) \leq x_1 \leq \phi(t,3), |x_2| \leq a_{t}(x_1), t< t_\varepsilon \}.
\end{equation}

On this domain, we have that
\begin{equation}\label{previousresult}
	v(t,\phi_\varepsilon(y)) = \chi_\varepsilon(y),
\end{equation}

where

\begin{equation}\label{phipd}
	\phi_\varepsilon(t,y) = y + t \frac{1+\chi_\varepsilon(y)}{1-\chi_\varepsilon(y)}.
\end{equation}

In the following, we show that $||v_\varepsilon(t)||_{H^{3/4}_{x_1}} \rightarrow \infty$ as $t \rightarrow t_\varepsilon$. 

Differentiating \ref{phipd} with respect to $y$ and $t$ we compute the derivatives that we will need,

\begin{equation}\label{phiderivatives}
	\begin{split}
		&\phi_{\varepsilon,y}(t,y) = 1+ 2 t \frac{\chi'_\varepsilon(y)}{(1-\chi_\varepsilon(y))^2}, \hspace{0.5cm} \phi_{\varepsilon,yy} =  2 t \frac{\chi''_\varepsilon(y)(1-\chi_\varepsilon(y)) +2 \chi'_\varepsilon(y)^2}{(1-\chi_\varepsilon(y))^3},\\
		&\phi_{\varepsilon,ty}(t,y) = 2 \frac{\chi'_\varepsilon(y)}{(1-\chi_\varepsilon(y))^2}, \hspace{0.5cm} \phi_{\varepsilon,tyy}(t,y) = 2 \frac{\chi''_\varepsilon(y)(1-\chi_\varepsilon(y)) + 2 \chi'_\varepsilon(y)^2}{(1-\chi_\varepsilon(y))^3}.\\
	\end{split}
\end{equation}

Now $y \mapsto \frac{|\chi'_\varepsilon(y)|}{(1-\chi_\varepsilon(y))^2}$ is a continuous function on a compact set, it reaches its maximum $M_\varepsilon$ at $y= \nu_\varepsilon$. Call $t_\varepsilon = \frac{1}{M_\varepsilon}$, we have the following properties

\begin{equation}
	\begin{split}
		&\phi_{\varepsilon,y}(t,y) \neq 0 \text{ for } t<t_\varepsilon, \hspace{0.5cm} \phi_{\varepsilon,y}(t_\varepsilon, \nu_\varepsilon) = 0, \hspace{0.5cm} t_\varepsilon \leq \frac{C}{|\ln(\varepsilon)|^{\alpha}}.
	\end{split}
\end{equation}

We have that for $y>\varepsilon$, $\phi_{\varepsilon,yy}>0$ which means that $\nu_\varepsilon < \varepsilon$. Note that $\phi_\varepsilon(t,\cdot)$ is an injection. We choose $\psi_\varepsilon$ such that $\nu_\varepsilon$ is unique and such that $\phi_{\varepsilon,yyy}(t_\varepsilon,\nu_\varepsilon) \neq 0$. Note that we hence have at $t=t_\varepsilon$ the following properties for $y$ close enough to $\nu_\varepsilon$ (we can choose the $\delta_\varepsilon$ such that this is satisfied since it only depends on $\varepsilon$ and not on $t$).

\begin{equation}\label{properties}
	\begin{split}
		&\phi_{\varepsilon,yy}(t_\varepsilon,\nu_\varepsilon) = 0,\\
		&\exists C_{1,\varepsilon},C_{2,\varepsilon} >0,\hspace{0.2cm} C_{1,\varepsilon} (y-\nu_\varepsilon) \leq \phi_{\varepsilon,yy}(t_\varepsilon,y) \leq C_{2,\varepsilon} (y-\nu_\varepsilon), \\
		&\exists C_{1,\varepsilon},C_{2,\varepsilon} >0,\hspace{0.2cm} C_{1,\varepsilon} (y-\nu_\varepsilon)^2 \leq \phi_{\varepsilon,y}(t_\varepsilon,y) \leq C_{2,\varepsilon} (y -\nu_\varepsilon)^2. 
	\end{split}
\end{equation} 

\begin{remark}
	The fact that $\phi_{\varepsilon,yy}$ is not of constant sign really is an issue for the estimation of the Sobolev norm, but we will explain later how we address this issue. $\phi_{\varepsilon,y}$ however, is of constant sign. 
\end{remark}

Also, we recall 

\begin{equation}\label{properties2}
	\exists M_\varepsilon<0, \hspace{0.2cm} 0 > \chi_\varepsilon(y) > M_\varepsilon, \hspace{0.2cm} 0 > \chi_\varepsilon'(y) > M_\varepsilon, \hspace{0.2cm} |\chi_\varepsilon''(y)| \leq |M_\varepsilon|
\end{equation}

Now, using Taylor expansions and the expressions given by (\ref{properties}) and (\ref{properties2}), we write the following inequalities

\begin{multline}\label{properties3}
	C_{1,\varepsilon} (\nu_\varepsilon-y)^2 + C_{1,\varepsilon} \phi_{\varepsilon,ty}(t_\varepsilon,\nu_\varepsilon)(t-t_\varepsilon) \\ \leq  \phi_{\varepsilon,y}(t,y) \leq C_{2,\varepsilon} (\nu_\varepsilon-y)^2 + C_{2,\varepsilon} \phi_{\varepsilon,ty}(t_\varepsilon,\nu_\varepsilon)(t-t_\varepsilon), \hspace{0.1cm} (i)\\
	C_{1,\varepsilon} (\nu_\varepsilon-y) \leq \phi_{\varepsilon,yy}(t,y) \leq C_{2,\varepsilon} (\nu_\varepsilon-y), \hspace{0.1cm} (ii)
\end{multline}

where $(ii)$ comes from the fact that $\phi_{\varepsilon,yy}(t,y) = \frac{t}{t_\varepsilon} \phi_{\varepsilon,yy}(t_\varepsilon,y)$. Now, using (\ref{properties}) and (\ref{properties3}), because $\phi_{\varepsilon,ty}<0$ and $(t-t_\varepsilon) <0$, both components have the same sign, and we can write

\begin{equation}
	C_{1,\varepsilon} |\nu_\varepsilon-y|^2 + C_{1,\varepsilon} (t_\varepsilon-t) \leq  |\phi_{\varepsilon,y}(t,y)| \leq C_{2,\varepsilon} |\nu_\varepsilon-y|^2 + C_{2,\varepsilon} (t_\varepsilon-t) \hspace{0.1cm},
\end{equation}

and $\phi_{\varepsilon,y}<0$ everywhere for $t<t_\varepsilon$.

Informally, we recall that our goal is to obtain lower and upper bounds for the integral (we omitted the independent variable)

\begin{equation}
	\int_{x_1} \int_{y} \frac{1}{|x_1-y|^{1/2-2 \lambda}} \frac{\partial^2{h_\varepsilon}}{\partial_{x_1}^2}(t,x_1) \frac{\partial^2{h_\varepsilon}}{\partial_{x_1}^2}(t,y),
\end{equation}

We obtain an estimation for the remaining terms.

\begin{equation}\label{remai}
	\begin{split}
		&h_\varepsilon(t,x_1,x_2) = \psi^1_{\varepsilon}(x_1) \psi_\varepsilon^2 (x_2) v_\varepsilon(t,x), \\
		&\frac{\partial h_\varepsilon}{\partial_{x_1}}(t,x_1,x_2) = \psi_\varepsilon^2(x_2) \big[ \psi_\varepsilon^{1,'} (x_1) v_\varepsilon(t,x_1) + \psi_\varepsilon^1(x_1) v_{\varepsilon,x_1}(t,x_1) \big],\\
		&\frac{\partial^2 h_\varepsilon}{\partial_{x_1}^2}(t,x_1,x_2) = \psi_{\varepsilon}^2(x_2) \big[ \psi_\varepsilon^{1,''}(x_1) v_\varepsilon(t,x_1) + 2 \psi_\varepsilon^{1,'}(x_1) v_{\varepsilon,x_1}(t,x_1) + \psi_\varepsilon^1(x_1) v_{\varepsilon,x_1x_1}(t,x_1) \big],\\  
		&\Rightarrow \exists C_{1,\varepsilon},C_{2,\varepsilon}>0, \hspace{0.2cm} C_{1,\varepsilon} \psi_\varepsilon^2(x_2) v_{\varepsilon,x_1x_1}(t,x_1) \leq \frac{\partial^2 h_\varepsilon}{\partial_{x_1}^2}(t,x_1,x_2) \leq C_{2,\varepsilon} \psi_\varepsilon^2(x_2) v_{\varepsilon,x_1x_1}(t,x_1) \hspace{0.1cm} (i) \\
		&\Rightarrow \exists C_{1,\varepsilon}>0, \hspace{0.2cm} \left|\frac{\partial^2 h_\varepsilon}{\partial_{x_1}^2}(t,x_1,x_2)\right| \leq C_{1,\varepsilon} \left|\psi_\varepsilon^2(x_2) v_{\varepsilon,x_1x_1}(t,x_1)\right|, \hspace{0.1cm} (ii) \\
	\end{split}
\end{equation}

where $(i)$ is valid when $x_1 \in [\phi_\varepsilon(t_\varepsilon,\nu_\varepsilon) - \delta_\varepsilon, \phi_\varepsilon(t_\varepsilon,\nu_\varepsilon) + \delta_\varepsilon]$ and $(ii)$ is valid everywhere. 

Now we are ready to start the proof of theorem \ref{conj}. We will use the formula obtained in lemma \ref{joachim} for the expression of the Sobolev norm.

\begin{proof}
	
	First, with $I_\varepsilon = \int_{x_2 = -2 \delta_\varepsilon}^{2 \delta_\varepsilon} \left( \psi_\varepsilon^2(x_2) \right)^2 dx_2$, we define the following three integrals.
	
	\begin{equation}
		\begin{split}
			I_\varepsilon(t)= \int_{x_1=\phi_\varepsilon(t_\varepsilon,\nu_\varepsilon) - 2 \delta_\varepsilon }^{\phi_\varepsilon(t_\varepsilon,\nu_\varepsilon) + 2 \delta_\varepsilon } &
			\int_{x_2=-2\delta_\varepsilon}^{2 \delta_\varepsilon} \Big[ \left( \frac{\partial^2 (v \psi_{\varepsilon}^1 \psi_\varepsilon^2)}{\partial x_1^2} \right)(t,x_1,x_2)  \\
			&\cdot \int_{y = \phi_\varepsilon(t_\varepsilon,\nu_\varepsilon) - 2 \delta_\varepsilon }^{\phi_\varepsilon(t_\varepsilon,\nu_\varepsilon) + 2 \delta_\varepsilon } |x_1-y|^{-1/2+2\lambda} \left( \frac{\partial^2 (v \psi_{\varepsilon}^1 \psi_\varepsilon^2)}{\partial x_1^2} \right)(t,y,x_2) dy \Big] dx_2 dx_1  \\
			= I_\varepsilon \int_{x_1=\phi_\varepsilon(t_\varepsilon,\nu_\varepsilon) - 2 \delta_\varepsilon }^{\phi_\varepsilon(t_\varepsilon,\nu_\varepsilon) - \delta_\varepsilon } & \Big[ \left( \frac{\partial^2 (v \psi_{\varepsilon}^1)}{\partial x_1^2} \right)(t,x_1)\\ &\cdot  \int_{y = \phi_\varepsilon(t_\varepsilon,\nu_\varepsilon) - 2 \delta_\varepsilon }^{\phi_\varepsilon(t_\varepsilon,\nu_\varepsilon) + 2 \delta_\varepsilon } |x_1-y|^{-1/2+2\lambda} \left( \frac{\partial^2 (v \psi_{\varepsilon}^1)}{\partial x_1^2} \right)(t,y) dy \Big] dx_1  \\
			+ I_\varepsilon \int_{x_1=\phi_\varepsilon(t_\varepsilon,\nu_\varepsilon) - \delta_\varepsilon }^{\phi_\varepsilon(t_\varepsilon,\nu_\varepsilon) + \delta_\varepsilon } & \Big[ \left( \frac{\partial^2 v}{\partial x_1^2} \right)(t,x_1) \\
			&\cdot  \int_{y = \phi_\varepsilon(t_\varepsilon,\nu_\varepsilon) - 2 \delta_\varepsilon }^{\phi_\varepsilon(t_\varepsilon,\nu_\varepsilon) + 2 \delta_\varepsilon } |x_1-y|^{-1/2+2\lambda} \left( \frac{\partial^2 v}{\partial x_1^2} \right)(t,y) dy \Big] dx_1  \\
			+ I_\varepsilon \int_{x_1=\phi_\varepsilon(t_\varepsilon,\nu_\varepsilon) + \delta_\varepsilon }^{\phi_\varepsilon(t_\varepsilon,\nu_\varepsilon) + 2 \delta_\varepsilon } & \Big[ \left( \frac{\partial^2 (v \psi_{\varepsilon}^1)}{\partial x_1^2} \right)(t,x_1) \\&\cdot  \int_{y = \phi_\varepsilon(t_\varepsilon,\nu_\varepsilon) - 2 \delta_\varepsilon }^{\phi_\varepsilon(t_\varepsilon,\nu_\varepsilon) + 2 \delta_\varepsilon }|x_1-y|^{-1/2+2\lambda} \left( \frac{\partial^2 (v \psi_{\varepsilon}^1)}{\partial x_1^2} \right)(t,y) dy \Big] dx_1  \\
			= I_\varepsilon (I_\varepsilon^1(t) + I_\varepsilon^2(t) +& I_\varepsilon^3(t))
		\end{split}
	\end{equation}
	
	\vspace{1cm}
	\begin{figure}[h]
		\begin{tikzpicture}
			\draw[semithick,->] (0,0,0) -- (0,0,3);
			\draw[semithick] (0,0,0) -- (0,0,-3);
			\draw (0,0,3.4) node {$x_2$};
			\draw[semithick,->] (0,0,0) -- (7,0,0);
			\draw (7.4,0,0) node {$x_1$};
			\draw[semithick,->] (0,0,0) -- (0,4,0);
			\draw (0,4.4,0) node {$t$};
			\draw (3,0,3.8) node {\footnotesize $x_1 = \nu_\varepsilon$};
			\draw[dashed] (3,0,-3) -- (3,0,3);
			\draw[dashed] (3,3,-3) -- (9,3,-3) -- (9,3,3) -- (3,3,3) -- (3,3,-3);
			\draw (10.2,3,-2) node {$P : t=t_\varepsilon$};
			\draw[dashed] (6,3,-3) -- (6,3,3);
			\draw (6,3,-3.4) node {\tiny $x_1= \phi_\varepsilon(t_\varepsilon,\nu_\varepsilon)$};
			\draw (7.5,3,-3.4) node {\tiny $+\delta_\varepsilon$};
			\draw (9,3,-3.4) node {\tiny $+2\delta_\varepsilon$};
			\draw (4.5,3,-3.4) node {\tiny $-\delta_\varepsilon$};
			\draw (3,3,-3.4) node {\tiny $-2\delta_\varepsilon$};
			\draw (2.6,3,3) node {\tiny $+2\delta_\varepsilon$};
			\draw (2.6,3,1.5) node {\tiny $+\delta_\varepsilon$};
			\draw (2.5,3,0) node {\tiny $x_2=0$};
			\draw (2.6,3,-1.5) node {\tiny $-\delta_\varepsilon$};
			
			\filldraw[green!15] (3.05,3,-2.95) -- (4.5,3,-2.95) -- (4.45,3,2.95) -- (3.05,3,2.95) --  (3.05,3,-2.95);\filldraw[orange!15] (7.5,3,-2.95) -- (9,3,-2.95) -- (9,3,2.95) -- (7.5,3,2.95) --  (7.5,3,-2.95);
			\filldraw[yellow!15] (4.5,3,-2.95) -- (5.95,3,-2.95) -- (5.95,3,2.95) -- (4.5,3,2.95) --  (4.5,3,-2.95);
			\filldraw[yellow!15] (6.05,3,-2.95) -- (6.05,3,2.95) -- (7.45,3,2.95) -- (7.45,3,-2.95) --  (6.05,3,-2.95);
			
			\draw[dashed,semithick,->] (3,0,0) parabola (6,2.8,0); 
			\draw (5.47,0.8,0) node {\tiny $x_1 = \phi_{\varepsilon}(t,\nu_\varepsilon)$};
			\draw[dashed] (4.5,3,-3) -- (4.5,3,3);
			\draw[dashed] (7.5,3,-3) -- (7.5,3,3);
			\draw[dashed] (3,3,-1.5) -- (9,3,-1.5);
			\draw[dashed] (3,3,1.5) -- (9,3,1.5);
			\draw[dashed] (3,3,0) -- (9,3,0);
			\draw (3.75,3,-2.25) node {\tiny $h_\varepsilon$};
			\draw (5.25,3,-2.25) node {\tiny $\psi_\varepsilon^1 = 1$};
			\draw (5.25,3,-0.75) node {\tiny $v_\varepsilon$};
			\draw (5.25,3,0.75) node {\tiny $v_\varepsilon$};
			\draw (5.25,3,2.25) node {\tiny $\psi_\varepsilon^1 = 1$};
			\draw (6.75,3,-2.25) node {\tiny $\psi_\varepsilon^1 = 1$};
			\draw (6.75,3,2.25) node {\tiny $\psi_\varepsilon^1 = 1$};
			
			\draw (6.75,3,-0.75) node {\tiny $v_\varepsilon$};
			\draw (6.75,3,0.75) node {\tiny $v_\varepsilon$};
			\draw (8.25,3,-2.25) node {\tiny $h_\varepsilon$};
			
			\draw (3.75,3,-0.75) node {\tiny $\psi_\varepsilon^2 = 1$};
			\draw (8.25,3,-0.75) node {\tiny $\psi_\varepsilon^2 = 1$};
			\draw (8.25,3,0.75) node {\tiny $\psi_\varepsilon^2 = 1$};
			\draw (3.75,3,0.75) node {\tiny $\psi_\varepsilon^2 = 1$};
			\draw (3.75,3,2.25) node {\tiny $h_\varepsilon$};
			\draw (8.25,3,2.25) node {\tiny $h_\varepsilon$};
			
			\filldraw[green!15] (8.5,2,0) rectangle (9,1.5,0);
			\draw (10.3,1.85,0) node {\tiny Domain of the first};
			\draw  (10.3,1.50,0) node {\tiny integral of $I_1$};
			\draw (10.3,0.85,0) node {\tiny Domain of the first};
			\draw  (10.3,0.50,0) node {\tiny integral of $I_2$};
			\draw (10.3,-0.15,0) node {\tiny Domain of the first};
			\draw  (10.3,-0.5,0) node {\tiny integral of $I_3$};
			
			\filldraw[yellow!15] (8.5,1,0) rectangle (9,0.5,0);
			
			\filldraw[orange!15] (8.5,0,0) rectangle (9,-0.5,0);
			\draw (6,-1,0) node {\footnotesize Note: The domain of $I_j(t)$ does not};
			\draw (6,-1.4,0) node {\footnotesize follow the line given by $\phi_\varepsilon$};
			\draw (6,-1.8,0) node {\footnotesize as $t$ varies, it "moves vertically". };
		\end{tikzpicture}
		\caption{Definition of $\psi_\varepsilon^1,\psi_\varepsilon^2,h_\varepsilon,I_\varepsilon^1,I_\varepsilon^2$ and $I_\varepsilon^3$.}
	\end{figure}
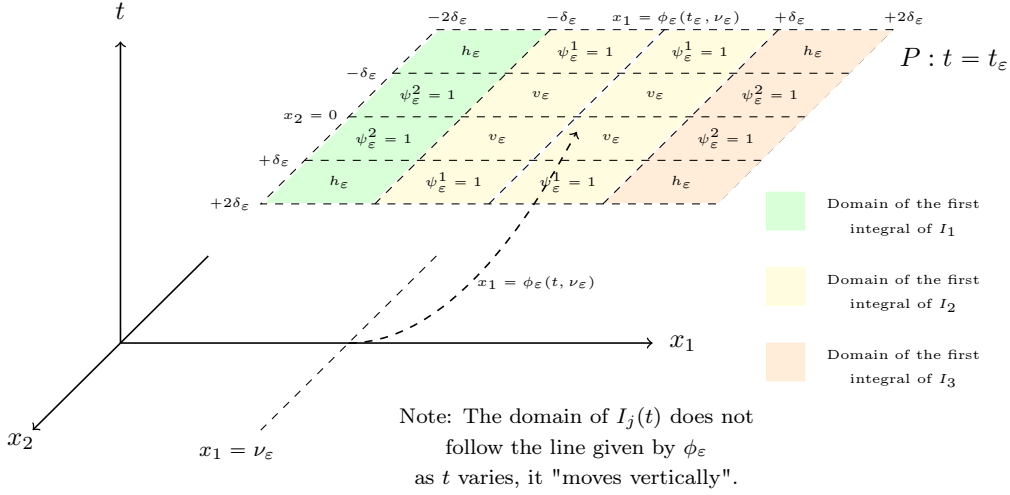
	
	The strategy of this proof will be to find upper bounds for $|I_\varepsilon^1|$ and $|I_\varepsilon^3|$ (as $t \rightarrow t_\varepsilon$) while finding a lower bound for $I_\varepsilon^2$ (as $t \rightarrow t_\varepsilon$) going to $+ \infty$ faster than the bounds on $|I_\varepsilon^1|$ and $|I_\varepsilon^3|$ are increasing, thus giving the convergence of $I_\varepsilon(t)$ to $+ \infty$ as $t\rightarrow t_\varepsilon$.
	
	Because of (\ref{remai}), we can study the integrals where $\frac{\partial^2 \left( v \psi_\varepsilon^1 \right)}{\partial_{x_1}^2}$ is replaced by $\frac{\partial^2 v}{\partial_{x_1}^2}$. Also, we make the change of variable (not explicitly relabelled) $x_1 = \phi_\varepsilon(t,x_1)$ and $y=\phi_\varepsilon(t,y)$, which makes the term $\phi_{\varepsilon,y}$ appear twice.
	
	Since
	
	\begin{equation}\label{v}
		v_\varepsilon(t,\phi_\varepsilon(t,y)) = \chi_\varepsilon(y),
	\end{equation}
	
	then,
	\begin{equation}\label{vx}
		\begin{split}
			v_{\varepsilon,x}(t,\phi_\varepsilon(t,y)) &= \frac{\chi'_{\varepsilon}(y)}{\phi_{\varepsilon,y}(t,y)} \\
			v_{\varepsilon,xx}(t,\phi_\varepsilon(t,y)) &= \frac{\chi''_{\varepsilon}(y)\phi_{\varepsilon,y}(t,y) - \chi'_\varepsilon(y) \phi_{\varepsilon,yy}(t,y) }{\phi^3_{\varepsilon,y}(t,y)}
		\end{split}
	\end{equation}

	Now by continuity, choose $\delta_\varepsilon$ so that for $z \in ] \phi_\varepsilon(t_\varepsilon,\nu_\varepsilon) - 2 \delta_\varepsilon, \phi_\varepsilon(t_\varepsilon,\nu_\varepsilon) + 2 \delta_\varepsilon [$, then $\phi_{\varepsilon,t}^{-1}(z) \in ] \varepsilon - \eta_\varepsilon, \varepsilon + \eta_\varepsilon [$, for a certain range of $t$ close enough to $t_\varepsilon$, of the form $]t_\varepsilon^0, t_\varepsilon[$.
	Hence, we can choose the $\delta_\varepsilon$ such that $\eta_\varepsilon$ is small enough, and all the results obtained in the preliminaries hold.
	
	For clarity, call
	
	\begin{equation}
		\begin{split}
			&\zeta_\varepsilon^1(t) = \phi_{\varepsilon,t}^{-1}(\phi_\varepsilon(t_\varepsilon,\nu_\varepsilon)-2\delta_\varepsilon)),\\
			&\zeta_\varepsilon^2(t) = \phi_{\varepsilon,t}^{-1}(\phi_\varepsilon(t_\varepsilon,\nu_\varepsilon)-\delta_\varepsilon)),\\
			&\zeta_\varepsilon^3(t) = \phi_{\varepsilon,t}^{-1}(\phi_\varepsilon(t_\varepsilon,\nu_\varepsilon)+\delta_\varepsilon)),\\
			&\zeta_\varepsilon^4(t) = \phi_{\varepsilon,t}^{-1}(\phi_\varepsilon(t_\varepsilon,\nu_\varepsilon)+2\delta_\varepsilon)),\\
		\end{split}
	\end{equation}
	
	and note that for $t_\varepsilon^1$ close enough to $t_\varepsilon$,
	
	\begin{equation}\label{order}
		\nu_\varepsilon- \eta_\varepsilon < \zeta_\varepsilon^1(t) < \zeta_\varepsilon^2(t) < \nu_\varepsilon < \zeta_\varepsilon^3(t) < \zeta_\varepsilon^4(t) < \nu_\varepsilon + \eta_\varepsilon.
	\end{equation}
	
	Also by continuity, for $t$ close enough to $t_\varepsilon$, there exist two constants $\zeta_\varepsilon^{2+}$ and $\zeta_\varepsilon^{3-}$ such that
	
	\begin{equation}\label{order5}
		\nu_\varepsilon- \eta_\varepsilon < \zeta_\varepsilon^1(t) < \zeta_\varepsilon^2(t) < \zeta_\varepsilon^{2+} < \nu_\varepsilon < \zeta_\varepsilon^{3-} < \zeta_\varepsilon^3(t) < \zeta_\varepsilon^4(t) < \nu_\varepsilon + \eta_\varepsilon.
	\end{equation}
	
	Let us now consider $I_\varepsilon^1$ (the case of $I_\varepsilon^3$ is similar). 
	
	\begin{multline}\label{I1lol}
		|I_\varepsilon^1(t)| \leq C_\varepsilon \int_{x_1=\zeta_\varepsilon^1(t)}^{\zeta_\varepsilon^2(t)} \int_{y=\zeta_\varepsilon^1(t)}^{\zeta_\varepsilon^4(t)} \frac{| \phi_{\varepsilon,y}(t,x_1) \phi_{\varepsilon,y}(t,y) |}{|\phi_\varepsilon(t,x_1) - \phi_\varepsilon(t,y)|^{1/2-2\lambda}} \left| v_{\varepsilon,xx}(t,x_1) v_{\varepsilon,xx}(t,y) \right| \\
		\leq C_\varepsilon \int_{x_1=\nu_\varepsilon-\eta_\varepsilon}^{\zeta_\varepsilon^{2+}} \left| \phi_{\varepsilon,y}(t,x_1) v_{\varepsilon,xx}(t,x_1) \right| \int_{y=\nu_\varepsilon - \eta_\varepsilon}^{\nu_\varepsilon+\eta_\varepsilon} \frac{|\phi_{\varepsilon,y}(t,y)|}{|\phi_\varepsilon(t,x_1) - \phi_\varepsilon(t,y)|^{1/2-2\lambda}} |v_{\varepsilon,xx}(t,y)|
	\end{multline}
	
	We study the inner integral of (\ref{I1lol}) for $x_1\in [\nu_\varepsilon - \eta_\varepsilon,\zeta_\varepsilon^{2+}]$.
	Now, both $\frac{1}{|\phi_\varepsilon(t,x_1) - \phi_\varepsilon(t,y)|^{1/2-2\lambda}}$ and $v_{\varepsilon,yy}(t,y)$ are unbounded in the second integral, but the regions where they are unbounded are uniformly disjoint in $t$ because of (\ref{order5}), so we split again the domain.
	
	\begin{multline}\label{I1lol2}
		\int_{y=\nu_\varepsilon - \eta_\varepsilon}^{\nu_\varepsilon+\eta_\varepsilon} \frac{1}{|\phi_\varepsilon(t,x_1) - \phi_\varepsilon(t,y)|^{1/2-2\lambda}} v_{\varepsilon,xx}(t,y) \leq_{(*1)} C_\varepsilon \int_{y=\nu_\varepsilon-\eta_\varepsilon}^{\zeta_\varepsilon^{2+}} \frac{1}{|\phi_\varepsilon(t,x_1) - \phi_\varepsilon(t,y)|^{1/2-2\lambda}} \\
		+ C_\varepsilon \int_{y=\zeta_\varepsilon^{2+}}^{\zeta_\varepsilon^{3-}} |\phi_{\varepsilon,y}(t,y) v_{\varepsilon,xx}(t,\phi_\varepsilon(t,y))| + C_\varepsilon \int_{y=\zeta_\varepsilon^{3-}}^{\nu_\varepsilon+\eta_\varepsilon} \frac{1}{|\phi_\varepsilon(t,x_1) - \phi_\varepsilon(t,y)|^{1/2-2\lambda}}\\
		\leq C_\varepsilon \left[ (i) + (ii) + (iii) \right]
	\end{multline}
	
	In $(*1)$, we used (\ref{properties}), (\ref{properties3}) and (\ref{vx}) and the fact that $|\nu_\varepsilon - y| \geq C_\varepsilon$ on the set $[\nu_\varepsilon - \eta_\varepsilon, \zeta_\varepsilon^{2+}]$.
	
	Now we give an upperbound for $(i)$ in (\ref{I1lol2}). The case $(iii)$ is trivial. 
	\begin{remark}
		For $I_\varepsilon^3(t)$, the case of $(i)$ is trivial, the $(ii)$ works the same as $(ii)$ for $I_\varepsilon^1(t)$, and $(iii)$ works the same as $(i)$ for $I_\varepsilon^1(t)$ that we do now.
	\end{remark}
	
	Now, using the mean value theorem as well as (\ref{properties3}), and $|c-\nu_\varepsilon| \geq |\zeta_\varepsilon^{2+} - \nu_\varepsilon|$, we obtain
	
	\begin{equation}
		(i) \leq C_\varepsilon \int_{y=\nu_\varepsilon - \eta_\varepsilon}^{\zeta_\varepsilon^{2+}} \frac{1}{|x_1-y|^{1/2-2\lambda}} \frac{1}{|\phi_{\varepsilon,y}(t,c)|^{1/2-2\lambda}} \leq C_\varepsilon \int_{y=\nu_\varepsilon - \eta_\varepsilon}^{\zeta_\varepsilon^{2+}} \frac{1}{|x_1-y|^{1/2-2\lambda}} \leq C_\varepsilon.
	\end{equation}
	
	Now, we study $(ii)$. Using the expression of $v_{\varepsilon,xx}$ given by (\ref{vx}), and the results obtained in the preliminary work, it is not clear if the bigger term is $\frac{1}{\phi_{\varepsilon,y}(t,x)^2}$ or $\frac{\phi_{\varepsilon,yy}(t,x)}{\phi_{\varepsilon,y}(t,x)^3}$. (The other factors being uniformly bounded for a fixed $\varepsilon$). Indeed, the root of $\phi_{\varepsilon,y}$ is of order $2$ in $x$ instead of the root being of order $1$ in $x$ for $\phi_{\varepsilon,yy}$, but the inequalities also involve a term depending on $t$ in $\phi_{\varepsilon,y}$. So we make the computations for both. 
	
	Using the expressions from (\ref{vx}) and (\ref{I1lol2}) as well as the inequalities obtained in (\ref{properties2}) and (\ref{properties3}), we obtain
	
	\begin{multline}\label{I1sol}
		(ii) \leq C_\varepsilon \int_{y=\zeta_\varepsilon^{2+}}^{\zeta_{\varepsilon}^{3-}} \frac{1}{\left| \phi_{\varepsilon,y}(t,y) \right|^2} + C_\varepsilon \int_{y=\zeta_\varepsilon^{2+}}^{\zeta_{\varepsilon}^{3-}} \frac{|\phi_{\varepsilon,yy}(t,y)|}{\left| \phi_{\varepsilon,y}(t,y) \right|^3} \leq \\
		C_\varepsilon \int_{y=\zeta_\varepsilon^{2+}}^{\zeta_{\varepsilon}^{3-}} \frac{1}{\left( |y-\nu_\varepsilon|^2 + (t_\varepsilon-t) \right)} + C_\varepsilon \int_{y=\zeta_\varepsilon^{2+}}^{\zeta_{\varepsilon}^{3-}} \frac{|y-\nu_\varepsilon|}{\left( |y-\nu_\varepsilon|^2 + (t_\varepsilon-t) \right)^2},
	\end{multline}

	and hence,
	
	\begin{equation}
		|(i) + (ii) + (iii)| \leq C_\varepsilon \frac{1}{(t_\varepsilon-t)}.
	\end{equation}
	
	Finally, we have the following upper bound for $I_\varepsilon^1(t)$, using $|x_1 - \nu_\varepsilon| \geq C_\varepsilon$,
	
	\begin{equation}\label{imaj}
		|I_\varepsilon^1(t)| \leq \frac{C_\varepsilon}{(t_\varepsilon-t)} \int_{x_1=\nu_\varepsilon - \eta_\varepsilon}^{\zeta_\varepsilon^{2-}} \left[ \frac{1}{ |\phi_{\varepsilon,y}(t,x_1)|} + \frac{|\phi_{\varepsilon,yy}(t,x_1)|}{|\phi_{\varepsilon,yy}(t,x_1)|^2} \right] \leq \frac{C_\varepsilon}{(t_\varepsilon-t)}.
	\end{equation}
	
	By symmetry, we also have $I_\varepsilon^3(t) \leq \frac{C_\varepsilon}{(t_\varepsilon-t)}$.
	
	Now, we do the estimation for $I_\varepsilon^2(t)$. 
	
	This one is more complicated because we do not have an upper bound for $\frac{\chi''(y)}{ (\phi_{\varepsilon,y}(t,x))^2} - \frac{\chi'(y) \phi_{\varepsilon,yy}(t,x)}{(\phi_{\varepsilon,y}(t,y))^3}$. Indeed, the first term is of constant sign whereas the second term is changing sign when $x$ is lower or bigger than $\nu_\varepsilon$. Also, the second term is not smaller than the first one. The idea we will use is the following one. In case of an integer Sobolev norm, the term is squared and is of constant sign. In the case of this fractional Sobolev norm, the kernel is not a Dirac but concentrates at $y=x_1$ as $\frac{1}{|\phi_{\varepsilon,y}(t,c)(x_1 - y)|^{1/2-2\lambda}}$.
	
	From (\ref{properties3}), we have 
	\begin{equation*}
		C_{1,\varepsilon} (\nu_\varepsilon - y) \leq \phi_{\varepsilon,yy}(t,y) \leq C_{2,\varepsilon} (\nu_\varepsilon - y).
	\end{equation*}
	
	We can write from (\ref{vx})
	
	\begin{multline}\label{integrand}
		\int_{x_1 = \zeta_\varepsilon^{2+}}^{\zeta_\varepsilon^{3-}} \int_{y = \zeta_\varepsilon^{2+}}^{\zeta_\varepsilon^{3-}} \frac{1}{|\phi_\varepsilon(t,x_1) - \phi_\varepsilon(t,y) |^{1/2-2\lambda}} \phi_{\varepsilon,y}(t,x_1) v_{xx}(t,x_1) \phi_{\varepsilon,y}(t,y) v_{xx}(t,y) \\
		= \int_{x_1 = \zeta_\varepsilon^{2+}}^{\zeta_\varepsilon^{3-}} \int_{y = \zeta_\varepsilon^{2+}}^{\zeta_\varepsilon^{3-}} \frac{1}{|\phi_\varepsilon(t,x_1) - \phi_\varepsilon(t,y) |^{1/2-2\lambda}} \Bigg[ \frac{\chi''(x_1)}{\phi_{\varepsilon,y}(t,x_1)} \frac{\chi''(y)}{\phi_{\varepsilon,y}(t,y)} - \frac{\chi''(x_1)}{\phi_{\varepsilon,y}(t,x_1)}\frac{\chi'(y) \phi_{\varepsilon,yy}(t,y)}{\phi_{\varepsilon,y}(t,y)^2} \\
		- \frac{\chi''(y)}{\phi_{\varepsilon,y}(t,y)}\frac{\chi'(x_1) \phi_{\varepsilon,yy}(t,x_1)}{\phi_{\varepsilon,y}(t,x_1)^2} + \frac{\chi'(x_1) \phi_{\varepsilon,yy}(t,x_1)}{\phi_{\varepsilon,y}(t,x_1)^2} \frac{\chi'(y) \phi_{\varepsilon,yy}(t,y)}{\phi_{\varepsilon,y}(t,y)^2}\Bigg] = (i) + (ii) + (iii) + (iv) 
	\end{multline}

	We will show that $(iv) >> |(i)| + |(ii)| + |(iii)|$ as $t\rightarrow t_\varepsilon$. Note that we do not study the integrals corresponding to $(x_1,y) \in [\zeta_{\varepsilon}^{2+}, \zeta_\varepsilon^{3-}] \times [\nu_\varepsilon-\eta_\varepsilon,\zeta_\varepsilon^{2+}]$ and $(x_1,y) \in [\zeta_{\varepsilon}^{2+}, \zeta_\varepsilon^{3-}] \times [\zeta_\varepsilon^{3-},\nu_\varepsilon+\eta_\varepsilon]$ because the method is identical to the one we used for $I_1$.
	
	First, we look at $(iv)$.
	In this case, we have to keep the two integrals together. Because $\int_{y} \frac{\chi'(y) \phi_{\varepsilon,yy}(t,y)}{\phi_{\varepsilon,y}(t,y)^2}$ is going to be smaller than $\int_{y} \frac{\chi''(y) }{\phi_{\varepsilon,y}(t,y)^1}$ because of $\phi_{yy}$ anti-symmetric properties (with $\nu_\varepsilon$ as the center). But the weight function will be higher when $(x-\nu_\varepsilon)$ and $(y-\nu_\varepsilon)$ are of the same sign. \\We denote $i(x_1,y) := \frac{\chi'(x_1) \phi_{\varepsilon,yy}(t,x_1)}{\phi_{\varepsilon,y}(t,x_1)^2} \frac{\chi'(y) \phi_{\varepsilon,yy}(t,y)}{\phi_{\varepsilon,y}(t,y)^2} \frac{1}{|\phi_\varepsilon(t,x_1) -\phi_\varepsilon(t,y)|^{1/2-2\lambda}}$. In figure \ref{alpha}, every pole of $i$ is displayed as a thick line. We see that the part corresponding to $\alpha$ and $\delta$ where $i\geq 0$ is going to be bigger than the part corresponding to $\beta$ and $\gamma$ where $ i \leq 0$. We also denote
	
	\begin{equation}\label{domains}
		\begin{split}
			&\alpha = \{(x_1,y) \in [\zeta_\varepsilon^{2+},\nu_\varepsilon] \times [\zeta_{\varepsilon}^{2+},\nu_\varepsilon] \}, \hspace{0.2cm} \beta = \{(x_1,y) \in [\nu_\varepsilon,\zeta_\varepsilon^{3-}] \times [\zeta_{\varepsilon}^{2+},\nu_\varepsilon] \}\\
			&\gamma = \{(x_1,y) \in [\zeta_\varepsilon^{2+},\nu_\varepsilon] \times [\nu_\varepsilon,\zeta_\varepsilon^{3-}] \}, \hspace{0.2cm} \delta = \{(x_1,y) \in [\nu_\varepsilon,\zeta_\varepsilon^{3-}] \times [\nu_\varepsilon,\zeta_\varepsilon^{3-}] \}
		\end{split}
	\end{equation}
	\begin{figure}[h]\centering 
		\begin{tikzpicture}[domain=0:6,scale=0.7] 
			\draw[very thin,color=gray] (-0.1,-1.1) grid (6.1,5.1);
			\draw[thick,->] (-0.2,2) -- (6.2,2) node[right] {$x_1 = \nu_\varepsilon$}; 
			\draw (-0.8,0.5) node {$\alpha, i \geq 0$};
			\draw (6.8,3.5) node {$\delta, i \geq 0$};
			\draw (6.8,0.5) node {$\beta, i \leq 0$};
			\draw (-0.8,3.5) node {$\gamma, i \leq 0$};
			\draw[thick,->] (3,-1.2) -- (3,5.2) node[above] {$y=\nu_\varepsilon$};
			\draw[thick]    plot (\x,\x-1)             node[right] {$y =x_1$}; 
		\end{tikzpicture}
		\caption{Definition of $\alpha$, $\beta$, $\gamma$ and $\delta$}\label{alpha}
	\end{figure}
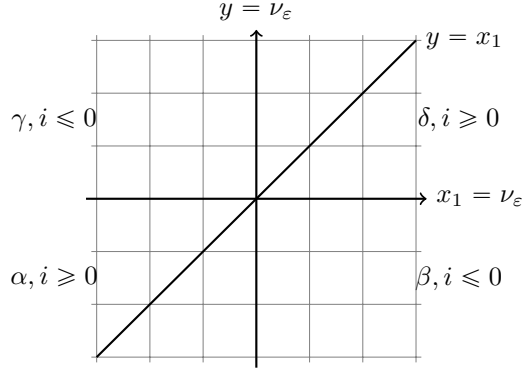
	
	And we have 
	\begin{equation}
		(iv) = \int\int_{\alpha} i(x_1,y) + \int\int_{\beta} i(x_1,y) +\int\int_{\gamma} i(x_1,y) + \int\int_{\delta} i(x_1,y). 
	\end{equation}
	
	We regroup the term $\delta$ and $\beta$ together. We will show that this term is non-negative, and also provide a lower bound for it. Without relabelling, we can choose $\zeta_\varepsilon^{2+}$ and $\zeta_\varepsilon^{3-}$ to be symmetric with respect to $\nu_\varepsilon$. We denote $\kappa_\varepsilon = \zeta_\varepsilon^{3-} - \nu_\varepsilon$.
	
	To do our estimate, we will need the following refined approximation for $\phi_{\varepsilon,yy}(t,y)$, that comes from a Taylor's expansion.
	
	\begin{equation}\label{refinement}
		C(\nu_\varepsilon-x_1) - C_1 (\nu_\varepsilon-x_1)^2 \leq \phi_{\varepsilon,yy}(t,x_1) \leq C(\nu_\varepsilon-x_1) + C_2 (\nu_\varepsilon - x_1)^2.
	\end{equation}
	
	\begin{multline}\label{Jdefinition}
		J=\int\int_{\delta \cup \beta}  \frac{1}{|\phi_\varepsilon(t,x_1) - \phi_\varepsilon(t,y) |^{1/2-2\lambda}} \frac{\chi'(x_1) \phi_{\varepsilon,yy}(t,x_1)}{\phi_{\varepsilon,y}(t,x_1)^2} \frac{\chi'(y) \phi_{\varepsilon,yy}(t,y)}{\phi_{\varepsilon,y}(t,y)^2} \\
		= \int_{x_1=0}^{\kappa_\varepsilon} \int_{y=0}^{\kappa_\varepsilon} \Bigg[ \frac{\chi'(\nu_\varepsilon+y) \phi_{\varepsilon,yy}(t,\nu_\varepsilon+y)\chi'(\nu_\varepsilon+x_1) \phi_{\varepsilon,yy}(t,\nu_\varepsilon+x_1)}{|\phi_\varepsilon(t,\nu_\varepsilon+x_1) - \phi_\varepsilon(t,\nu_\varepsilon+y) |^{1/2-2\lambda}\phi_{\varepsilon,y}(t,\nu_\varepsilon+x_1)^2\phi_{\varepsilon,y}(t,\nu_\varepsilon+y)^2}  \\
		+\frac{\chi'(\nu_\varepsilon-y) \phi_{\varepsilon,yy}(t,\nu_\varepsilon-y)\chi'(\nu_\varepsilon+x_1) \phi_{\varepsilon,yy}(t,\nu_\varepsilon+x_1)}{|\phi_\varepsilon(t,\nu_\varepsilon+x_1) - \phi_\varepsilon(t,\nu_\varepsilon-y) |^{1/2-2\lambda}\phi_{\varepsilon,y}(t,\nu_\varepsilon-y)^2\phi_{\varepsilon,y}(t,\nu_\varepsilon+x_1)^2}
		\Bigg]
	\end{multline}
	
	Now, we write $J$ defined in (\ref{Jdefinition}) as a difference term, using the approximation (\ref{refinement}).
	
	\begin{equation}
		J= J_1 + J_2 + J_3 + J_4,
	\end{equation}
	
	There exists $f$ a bounded function such that
	
	\begin{multline}
		J_1= C \int_{x_1=0}^{\kappa_\varepsilon} \int_{y=0}^{\kappa_\varepsilon} \Bigg[ \frac{\chi'(\nu_\varepsilon+x_1) \cdot x_1 \cdot \chi'(\nu_\varepsilon+y) \cdot y}{|\phi_\varepsilon(t,\nu_\varepsilon+x_1) - \phi_\varepsilon(t,\nu_\varepsilon+y) |^{1/2-2\lambda}\phi_{\varepsilon,y}(t,\nu_\varepsilon+y)^2\phi_{\varepsilon,y}(t,\nu_\varepsilon+x_1)^2}\\
		-\frac{\chi'(\nu_\varepsilon-y) \cdot y \cdot \chi'(\nu_\varepsilon+x_1) \cdot x_1}{|\phi_\varepsilon(t,\nu_\varepsilon+x_1) - \phi_\varepsilon(t,\nu_\varepsilon-y) |^{1/2-2\lambda}\phi_{\varepsilon,y}(t,\nu_\varepsilon-y)^2\phi_{\varepsilon,y}(t,\nu_\varepsilon+x_1)^2} \Bigg], \\
	\end{multline}

	\begin{multline}
		J_2= C \int_{x_1=0}^{\kappa_\varepsilon} \int_{y=0}^{\kappa_\varepsilon} \Bigg[ \frac{\chi'(\nu_\varepsilon+y) \cdot y\cdot \chi'(\nu_\varepsilon+x_1) \cdot f(\nu_\varepsilon+x_1) \cdot x_1^2}{|\phi_\varepsilon(t,\nu_\varepsilon+x_1) - \phi_\varepsilon(t,\nu_\varepsilon+y) |^{1/2-2\lambda}\phi_{\varepsilon,y}(t,\nu_\varepsilon+y)^2\phi_{\varepsilon,y}(t,\nu_\varepsilon+x_1)^2} \\
		-\frac{\chi'(\nu_\varepsilon-y) \cdot  y \cdot \chi'(\nu_\varepsilon+x_1) \cdot f(\nu_\varepsilon+x_1) \cdot x_1^2}{|\phi_\varepsilon(t,\nu_\varepsilon+x_1) - \phi_\varepsilon(t,\nu_\varepsilon-y) |^{1/2-2\lambda}\phi_{\varepsilon,y}(t,\nu_\varepsilon-y)^2\phi_{\varepsilon,y}(t,\nu_\varepsilon+x_1)^2}\Bigg], \\
		J_3= C \int_{x_1=0}^{\kappa_\varepsilon} \int_{y=0}^{\kappa_\varepsilon} \Bigg[ \frac{\chi'(\nu_\varepsilon+y) \cdot f(\nu_\varepsilon+y) \cdot y^2 \cdot \chi'(\nu_\varepsilon+x_1) \cdot x_1}{|\phi_\varepsilon(t,\nu_\varepsilon+x_1) - \phi_\varepsilon(t,\nu_\varepsilon+y) |^{1/2-2\lambda}\phi_{\varepsilon,y}(t,\nu_\varepsilon+y)^2\phi_{\varepsilon,y}(t,\nu_\varepsilon+x_1)^2}\\
		-\frac{\chi'(\nu_\varepsilon-y) \cdot f(\nu_\varepsilon+y) \cdot y^2 \cdot \chi'(\nu_\varepsilon+x_1) \cdot x_1}{|\phi_\varepsilon(t,\nu_\varepsilon+x_1) - \phi_\varepsilon(t,\nu_\varepsilon-y) |^{1/2-2\lambda}\phi_{\varepsilon,y}(t,\nu_\varepsilon-y)^2\phi_{\varepsilon,y}(t,\nu_\varepsilon+x_1)^2}\Bigg], \\
		J_4= C \int_{x_1=0}^{\kappa_\varepsilon} \int_{y=0}^{\kappa_\varepsilon} \Bigg[ \frac{\chi'(\nu_\varepsilon+y) \cdot f(\nu_\varepsilon+y) \cdot y^2 \cdot \chi'( \nu_\varepsilon+ x_1) \cdot f(\nu_\varepsilon+x_1)\cdot x_1^2}{|\phi_\varepsilon(t,\nu_\varepsilon+x_1) - \phi_\varepsilon(t,\nu_\varepsilon+y) |^{1/2-2\lambda}\phi_{\varepsilon,y}(t,\nu_\varepsilon+y)^2\phi_{\varepsilon,y}(t,\nu_\varepsilon+x_1)^2} \\
		-\frac{\chi'(\nu_\varepsilon-y) \cdot f(\nu_\varepsilon+y) \cdot y^2 \cdot \chi'(\nu_\varepsilon+x_1) \cdot f(\nu_\varepsilon+x_1) \cdot x_1^2}{|\phi_\varepsilon(t,\nu_\varepsilon+x_1) - \phi_\varepsilon(t,\nu_\varepsilon-y) |^{1/2-2\lambda}\phi_{\varepsilon,y}(t,\nu_\varepsilon-y)^2\phi_{\varepsilon,y}(t,\nu_\varepsilon+x_1)^2} \Bigg]. \\
	\end{multline}
	
	The main contribution will come from $J_1$. We will show later that the contributions from $J_2$, $J_3$ and $J_4$ are smaller. Let us, for now, focus only on $J_1$. We decompose the integrand as follows.
	
	\begin{multline}
		\frac{1}{|\phi_\varepsilon(t,\nu_\varepsilon+x_1) - \phi_\varepsilon(t,\nu_\varepsilon+y) |^{1/2-2\lambda}} \frac{\chi'(\nu_\varepsilon+x_1) \cdot x_1}{\phi_{\varepsilon,y}(t,\nu_\varepsilon+x_1)^2} \frac{\chi'(\nu_\varepsilon+y) \cdot y}{\phi_{\varepsilon,y}(t,\nu_\varepsilon+y)^2} \\
		-\frac{1}{|\phi_\varepsilon(t,\nu_\varepsilon+x_1) - \phi_\varepsilon(t,\nu_\varepsilon-y) |^{1/2-2\lambda}} \frac{\chi'(\nu_\varepsilon+x_1) \cdot x_1}{\phi_{\varepsilon,y}(t,\nu_\varepsilon+x_1)^2} \frac{\chi'(\nu_\varepsilon-y) \cdot y}{\phi_{\varepsilon,y}(t,\nu_\varepsilon-y)^2} \\
		= \Bigg[ \frac{1}{|\phi_\varepsilon(t,\nu_\varepsilon+x_1) - \phi_\varepsilon(t,\nu_\varepsilon+y) |^{1/2-2\lambda}} \frac{\chi'(\nu_\varepsilon+x_1) \cdot x_1}{\phi_{\varepsilon,y}(t,\nu_\varepsilon+x_1)^2} \frac{\chi'(\nu_\varepsilon+y) \cdot y}{\phi_{\varepsilon,y}(t,\nu_\varepsilon+y)^2} \\
		- \frac{1}{|\phi_\varepsilon(t,\nu_\varepsilon+x_1) - \phi_\varepsilon(t,\nu_\varepsilon-y) |^{1/2-2\lambda}} \frac{\chi'(\nu_\varepsilon+x_1) \cdot x_1}{\phi_{\varepsilon,y}(t,\nu_\varepsilon+x_1)^2} \frac{\chi'(\nu_\varepsilon+y) \cdot y}{\phi_{\varepsilon,y}(t,\nu_\varepsilon+y)^2} \Bigg] \\
		+ \Bigg[ \frac{1}{|\phi_\varepsilon(t,\nu_\varepsilon+x_1) - \phi_\varepsilon(t,\nu_\varepsilon-y) |^{1/2-2\lambda}} \frac{\chi'(\nu_\varepsilon+x_1) \cdot x_1}{\phi_{\varepsilon,y}(t,\nu_\varepsilon+x_1)^2} \frac{\chi'(\nu_\varepsilon+y) \cdot y}{\phi_{\varepsilon,y}(t,\nu_\varepsilon+y)^2} \\
		- \frac{1}{|\phi_\varepsilon(t,\nu_\varepsilon+x_1) - \phi_\varepsilon(t,\nu_\varepsilon-y) |^{1/2-2\lambda}} \frac{\chi'(\nu_\varepsilon+x_1) \cdot x_1}{\phi_{\varepsilon,y}(t,\nu_\varepsilon+x_1)^2} \frac{\chi'(\nu_\varepsilon+y) \cdot y}{\phi_{\varepsilon,y}(t,\nu_\varepsilon-y)^2} \Bigg] \\
		+ \Bigg[ \frac{1}{|\phi_\varepsilon(t,\nu_\varepsilon+x_1) - \phi_\varepsilon(t,\nu_\varepsilon-y) |^{1/2-2\lambda}} \frac{\chi'(\nu_\varepsilon+x_1) \cdot x_1}{\phi_{\varepsilon,y}(t,\nu_\varepsilon+x_1)^2} \frac{\chi'(\nu_\varepsilon+y) \cdot y}{\phi_{\varepsilon,y}(t,\nu_\varepsilon-y)^2} \\
		- \frac{1}{|\phi_\varepsilon(t,\nu_\varepsilon+x_1) - \phi_\varepsilon(t,\nu_\varepsilon-y) |^{1/2-2\lambda}} \frac{\chi'(\nu_\varepsilon+x_1) \cdot x_1}{\phi_{\varepsilon,y}(t,\nu_\varepsilon+x_1)^2} \frac{\chi'(\nu_\varepsilon-y) \cdot y}{\phi_{\varepsilon,y}(t,\nu_\varepsilon-y)^2} \Bigg] \\
		=D_1 + D_2 + D_3.
	\end{multline}
	
	To shorten a bit the notations, we will denote 
	
	\begin{equation}
		\alpha_{\pm} = \left| \phi_\varepsilon(t,\nu_\varepsilon+x_1) - \phi_\varepsilon(t,\nu_\varepsilon\pm y) \right|.
	\end{equation}
	
	For $D_1$, we write
	
	\begin{multline}
		\frac{1}{|\phi_\varepsilon(t,\nu_\varepsilon+x_1) - \phi_\varepsilon(t,\nu_\varepsilon+y) |^{1/2-2\lambda}} - \frac{1}{|\phi_\varepsilon(t,\nu_\varepsilon+x_1) - \phi_\varepsilon(t,\nu_\varepsilon-y) |^{1/2-2\lambda}} \\
		= \frac{\alpha_- - \alpha_+}{\alpha_+^{1/2-2\lambda} \alpha_-^{1/2-2\lambda}\left(\alpha_+^{1/2+2\lambda} + \alpha_-^{1/2+2\lambda}\right)} + \frac{\alpha_+^{4\lambda} - \alpha_-^{4\lambda}}{\left(\alpha_+^{1/2+2\lambda} + \alpha_-^{1/2+2\lambda}\right)}
	\end{multline}
	
	It is clear that $D_1$ is nonnegative when $x_1\geq y$. We now provide a lower bound for $D_{1,x_1 \leq y}$. We will show that it is positive up to a smaller order term, and provide a lower bound for the positive term. We hence now consider $x_1 \leq y$.
	
	We write, 
	
	\begin{multline}\label{lowerlol}
		\alpha_- - \alpha_+ = \left| \phi_\varepsilon(t,\nu_\varepsilon+x_1)-\phi_\varepsilon(t,\nu_\varepsilon-y) \right| - \left| \phi_\varepsilon(t,\nu_\varepsilon+x_1) - \phi_\varepsilon(t,\nu_\varepsilon+y) \right| \\
		= 2 \phi_\varepsilon(t,\nu_\varepsilon+x_1) - \phi_\varepsilon(t,\nu_\varepsilon+y) - \phi_\varepsilon(t,\nu_\varepsilon-y) = \int_{s=-y}^{x_1} \phi_{\varepsilon,y}(t,\nu_\varepsilon+s) - \int_{s=x_1}^y \phi_{\varepsilon,y}(t,\nu_\varepsilon+s) \\
		= \int_{s=-x_1}^{x_1} \phi_{\varepsilon,y}(t,\nu_\varepsilon+s) + \int_{s=x_1}^y \left( \phi_{\varepsilon,y} (t,\nu_\varepsilon-s) - \phi_{\varepsilon,y}(t,\nu_\varepsilon+s) \right).
	\end{multline}
	
	We will use this idea to provide a lower bound for $D_1$. 
	
	\begin{multline}\label{unlabellabel}
		\frac{1}{\alpha_+^{1/2-2\lambda}} - \frac{1}{\alpha_-^{1/2-2\lambda}} = \int_{s=\alpha_-}^{\alpha_+} -(\frac{1}{2}-2\lambda) \frac{1}{s^{3/2-2\lambda}} = C \int_{s=\alpha_+}^{\alpha_-} s^{-3/2+2\lambda}\\
		= C \int_{s=\alpha_+}^{\alpha_+ + \int_{z=-x_1}^{x_1} \phi_{\varepsilon,y}(t,\nu_\varepsilon+z) } s^{-3/2+2\lambda} \\
		+ C\int_{s=\alpha_+ + \int_{z=-x_1}^{x_1}\phi_{\varepsilon,y}(t,\nu_\varepsilon+z) }^{\alpha_+ + \int_{z=-x_1}^{x_1} \phi_{\varepsilon,y}(t,\nu_\varepsilon+z) +\int_{z=x_1}^y \left( \phi_{\varepsilon,y} (t,\nu_\varepsilon-z) - \phi_{\varepsilon,y}(t,\nu_\varepsilon+z) \right)} s^{-3/2+2\lambda} = (i) + (ii),
	\end{multline}
	
	where $(i)$ is nonnegative and $(ii)$ is small. First, we make the following upper bound for $|(ii)|$. If $\int_{x_1}^y  \left( \phi_{\varepsilon,y}(t,\nu_\varepsilon-s) - \phi_{\varepsilon,y}(t,\nu_\varepsilon+s) \right)>0$, then $(ii)$ is nonnegative. Otherwise, we have $\alpha_{+}+\int_{-x_1}^{x_1} \phi_{\varepsilon,y}(t,\nu_\varepsilon+s)+\int_{x_1}^y  \left( \phi_{\varepsilon,y}(t,\nu_\varepsilon-s) - \phi_{\varepsilon,y}(t,\nu_\varepsilon+s) \right)\leq \alpha_{+}+\int_{-x_1}^{x_1} \phi_{\varepsilon,y}(t,\nu_\varepsilon+s)$. We then proceed as follows 
	
	\begin{multline}
		\left| \int_{s=\alpha_{+}+\int_{-x_1}^{x_1} \phi_{\varepsilon,y}(t,\nu_\varepsilon+s)}^{\alpha_{+}+\int_{z=-x_1}^{x_1} \phi_{\varepsilon,y}(t,\nu_\varepsilon+z)+\int_{z=x_1}^y  \left( \phi_{\varepsilon,y}(t,\nu_\varepsilon-z) - \phi_{\varepsilon,y}(t,\nu_\varepsilon+z) \right)}  s^{-3/2-2\lambda} \right| \\
		= \left| \int_{s=\alpha_-}^{\alpha_- - \int_{z=x_1}^y  \left( \phi_{\varepsilon,y}(t,\nu_\varepsilon-z) - \phi_{\varepsilon,y}(t,\nu_\varepsilon+z) \right)} s^{-3/2+2\lambda} \right| \\
		\leq \left| \int_{z=x_1}^y  \left( \phi_{\varepsilon,y}(t,\nu_\varepsilon-z) - \phi_{\varepsilon,y}(t,\nu_\varepsilon+z) \right) \right| \cdot \alpha_-^{-3/2+2\lambda}. 
	\end{multline}
	
	From the Taylor expansion, we obtain for $x,y$ small enough (only depending on $\varepsilon$),
	
	\begin{equation}
		\left| \int_{z=x_1}^y  \left( \phi_{\varepsilon,y}(t,\nu_\varepsilon-s) - \phi_{\varepsilon,y}(t,\nu_\varepsilon+s) \right) \right| \leq C_1 \left( y^4 -x_1^4 + (t_\varepsilon-t)^2 + (t-t_\varepsilon)(y^2-x_1^2) \right). 
	\end{equation}
	
	This means that we obtain (up to a nonnegative contribution)
	
	\begin{equation}
		|(ii)| \leq C \left(y^4 + x_1^4 + (t_\varepsilon-t)\left( x_1^2 + y^2 \right) + (t_\varepsilon-t)^2 \right) \cdot \frac{1}{\alpha_-^{3/2-2\lambda}}.
	\end{equation}

	We now go on with the lower bound for $(i)$. Using the expression of $(i)$ provided by (\ref{unlabellabel}) as well as the mean value theorem, we obtain
	
	\begin{multline}\label{voiciunlabelnote}
		(i) = \int_{s=\alpha_+}^{\alpha_+ + \int_{z=-x_1}^{x_1} \phi_{\varepsilon,y}(t,\nu_\varepsilon+z) } s^{-3/2+2\lambda} \\
		\geq \left( \int_{z=-x_1}^{x_1} \phi_{\varepsilon,y}(t,\nu_\varepsilon+z) \right) \cdot \frac{1}{ \left(\phi_\varepsilon(t,\nu_\varepsilon+y) - \phi_{\varepsilon} (t,\nu_\varepsilon-x_1) \right)^{3/2-2\lambda}}. \\
		\geq C x_1 \phi_{\varepsilon,y}(t,c_1) \cdot \frac{1}{(y+x_1)^{3/2-2\lambda} \cdot \phi_{\varepsilon,y}(t,c_2)^{3/2-2\lambda}}.
	\end{multline}
	
	Now, since we have $x_1\leq y$, we obtain from (\ref{properties3})

	\begin{equation}\label{voiciunlabelnote1}
		\phi_{\varepsilon,y}(t,c_1) \geq C (t_\varepsilon-t),
	\end{equation}
	
	and
	
	\begin{equation}\label{voiciunlabelnote2}
		\phi_{\varepsilon,y}(t,c_2) \leq C \left( (t_\varepsilon-t) + y^2\right).
	\end{equation}
	
	Using (\ref{voiciunlabelnote1}) and (\ref{voiciunlabelnote2}) inside of (\ref{voiciunlabelnote}), we obtain
	
	\begin{equation}
		(i) \geq C \frac{x_1 (t_\varepsilon-t)}{(y+x_1)^{3/2-2\lambda} \cdot ((t_\varepsilon-t) + y^2)^{3/2-2\lambda} }.
	\end{equation}

	Now, we obtain for $D_1$,
	
	\begin{multline}\label{integrand2}
		D_1 = D_{x_1\geq y} + D_{x_1 \leq y} \geq D_{x_1\leq y} \\
		\geq C \int\int_{x_1\leq y} \frac{x_1 (t_\varepsilon-t)}{(y+x_1)^{3/2-2\lambda} \cdot ((t_\varepsilon-t) + y^2)^{3/2-2\lambda} } \frac{\chi'(\nu_\varepsilon+x_1) \cdot x_1}{\phi_{\varepsilon,y}(t,\nu_\varepsilon+x_1)^2} \frac{\chi'(\nu_\varepsilon+y) \cdot y}{\phi_{\varepsilon,y}(t,\nu_\varepsilon+y)^2} \\
		- C \int_{x_1=0}^{\kappa_\varepsilon}\int_{y=x_1}^{\kappa_\varepsilon} \frac{x^4}{\alpha_{-}^{3/2-2\lambda}}\frac{\chi'(\nu_\varepsilon+x_1) \cdot x_1}{\phi_{\varepsilon,y}(t,\nu_\varepsilon+x_1)^2} \frac{\chi'(\nu_\varepsilon+y) \cdot y}{\phi_{\varepsilon,y}(t,\nu_\varepsilon+y)^2} \\
		- C \int_{x_1=0}^{\kappa_\varepsilon}\int_{y=x_1}^{\kappa_\varepsilon} \frac{y^4}{\alpha_{-}^{3/2-2\lambda}}\frac{\chi'(\nu_\varepsilon+x_1) \cdot x_1}{\phi_{\varepsilon,y}(t,\nu_\varepsilon+x_1)^2} \frac{\chi'(\nu_\varepsilon+y) \cdot y}{\phi_{\varepsilon,y}(t,\nu_\varepsilon+y)^2} \\
		- C \int_{x_1=0}^{\kappa_\varepsilon} \int_{y=x_1}^{\kappa_\varepsilon} \frac{(t_\varepsilon-t)\cdot x_1^2}{\alpha_{-}^{3/2-2\lambda}}\frac{\chi'(\nu_\varepsilon+x_1) \cdot x_1}{\phi_{\varepsilon,y}(t,\nu_\varepsilon+x_1)^2} \frac{\chi'(\nu_\varepsilon+y) \cdot y}{\phi_{\varepsilon,y}(t,\nu_\varepsilon+y)^2} \\
		- C \int_{x_1=0}^{\kappa_\varepsilon} \int_{y=x_1}^{\kappa_\varepsilon} \frac{(t_\varepsilon-t)\cdot y^2}{\alpha_{-}^{3/2-2\lambda}}\frac{\chi'(\nu_\varepsilon+x_1) \cdot x_1}{\phi_{\varepsilon,y}(t,\nu_\varepsilon+x_1)^2} \frac{\chi'(\nu_\varepsilon+y) \cdot y}{\phi_{\varepsilon,y}(t,\nu_\varepsilon+y)^2} \\
		- C \int_{x_1=0}^{\kappa_\varepsilon} \int_{y=x_1}^{\kappa_\varepsilon} \frac{(t_\varepsilon-t)^2}{\alpha_{-}^{3/2-2\lambda}}\frac{\chi'(\nu_\varepsilon+x_1) \cdot x_1}{\phi_{\varepsilon,y}(t,\nu_\varepsilon+x_1)^2} \frac{\chi'(\nu_\varepsilon+y) \cdot y}{\phi_{\varepsilon,y}(t,\nu_\varepsilon+y)^2} \\
		= A_1 - B_1 - B_2 - B_3 - B_4 - B_5.
	\end{multline}
	
	Note that the proof also works for $\lambda=0$. We will now show that $A_1 \rightarrow \infty$, and that $J_2$, $J_3$, $J_4$, $D_2$, $D_3$, $B_1$, $B_2$, $B_3$, $B_4$ and $B_5$ are of a smaller order. We first consider $A_1$.
	
	Using $|\phi_{\varepsilon,y}(t,c(x_1,y))| \leq M_\varepsilon \left((x_1)^2 + (t_\varepsilon-t)\right)$, on $y > x_1$, we have from (\ref{integrand2}) and the new change of variable $(x,y) = (r-z,r+z)$, 
	
	\begin{multline}\label{integrand3}
		A_1 \geq  C_\varepsilon \int_{r=0}^{\kappa_\varepsilon/2} \int_{z=0}^r \frac{(r-z)(t_\varepsilon-t)}{r^{3/2-2\lambda}} \frac{r+z}{\left( \left( r+z \right)^2 + (t_\varepsilon-t) \right)^{7/2-2\lambda}} \frac{r-z}{\left( \left( r-z \right)^2 + (t_\varepsilon-t) \right)^{2}} \\
		\geq C_\varepsilon \int_{r=0}^{\kappa_\varepsilon/2} \frac{1}{r^{1/2-2\lambda}} \frac{(t_\varepsilon-t)}{\left( \left( 2r \right)^2 + (t_\varepsilon-t) \right)^{11/2-2\lambda}} \int_{z=0}^r (r-z)^2 \\
		\geq C_\varepsilon \int_{r=0}^{\kappa_\varepsilon/2} \frac{(t_\varepsilon-t)\cdot r^{5/2+2\lambda}}{\left( \left( 2r \right)^2 + (t_\varepsilon-t) \right)^{11/2-2\lambda}} \geq \frac{C_\varepsilon \cdot (t_\varepsilon-t)}{(t_\varepsilon-t)^{15/4-3\lambda}}\\
		\geq \frac{C_\varepsilon}{(t_\varepsilon-t)^{11/4-3\lambda}}.
	\end{multline}
	
	We now proceed with $B_1$. The case of $B_2$ is similar. We have since $\phi_{\varepsilon,y}(t,x) \geq (t_\varepsilon-t)$, with the mean value theorem
	
	\begin{multline}\label{lacestunb1lol}
		|B_1| \leq \frac{C_\varepsilon}{(t_\varepsilon-t)^{3/2-2\lambda}} \int_{r=0}^{\kappa_\varepsilon/2} \int_{z=0}^r  \frac{(r-z)^4}{r^{3/2-2\lambda} } \\
		\cdot  \frac{r+z}{\left( \left( r+z \right)^2 + (t_\varepsilon-t) \right)^{2}} \frac{r-z}{\left( \left( r-z \right)^2 + (t_\varepsilon-t) \right)^{2}} \\
		\leq \frac{C_\varepsilon}{(t_\varepsilon-t)^{3/2-2\lambda}} \int_{r=0}^{\kappa_\varepsilon/2}  \frac{r^4}{r^{1/2-2\lambda} \left(r^2+(t_\varepsilon-t)\right)^2} \int_{z=0}^r \frac{r-z}{ \left( (r-z)^2 + (t_\varepsilon-t) \right)^2}.
	\end{multline}

	Now, because
	
	\begin{equation}
		\int_{s=0}^\infty \frac{s}{\left(s^2 + (t_\varepsilon-t) \right)^2} \leq \frac{C}{(t_\varepsilon-t)},
	\end{equation}
	
	(\ref{lacestunb1lol}) yields
	
	\begin{equation}
		|B_1| \leq \frac{C_\varepsilon}{(t_\varepsilon-t)^{5/2-2\lambda}} \int_{r=0}^{\kappa_\varepsilon/2}  \frac{r^{7/2+2\lambda}}{ \left(r^2+(t_\varepsilon-t)\right)^2} \leq \frac{C_\varepsilon}{(t_\varepsilon-t)^{10/4-2\lambda}}.
	\end{equation}
	
	The cases of $B_3$, $B_4$ and $B_5$ are similarly done.
		
	This concludes the proof for the first difference term $D_1$. Overall, we obtain
	
	\begin{equation}
		D_1 \geq \frac{C_\varepsilon}{(t_\varepsilon-t)^{11/4-3\lambda}}.
	\end{equation}
	
	Now, we study $D_2$ and $D_3$.
	
	For $D_2$, we first study the term
	
	\begin{multline}
		\frac{1}{\phi_{\varepsilon,y}(t,\nu_\varepsilon+y)^2} - \frac{1}{\phi_{\varepsilon,y}(t,\nu_\varepsilon-y)^2} \\
		= \frac{\left(\phi_{\varepsilon,y}(t,\nu_\varepsilon+y) + \phi_{\varepsilon,y}(t,\nu_\varepsilon-y) \right) \cdot \left( \phi_{\varepsilon,y}(t,\nu_\varepsilon-y) - \phi_{\varepsilon,y}(t,\nu_\varepsilon+y) \right)}{\phi_{\varepsilon,y}(t,\nu_\varepsilon+y)^2 \phi_y(t,\nu_\varepsilon-y)^2}.
	\end{multline}
	
	A Taylor expansion of $\phi_{\varepsilon,y}$ gives that $\phi_{\varepsilon,y}(t,\nu_\varepsilon-y) - \phi_{\varepsilon,y}(t,\nu_\varepsilon+y) = f_1(t,y) y^3 + f_2(t,y) (t-t_\varepsilon)^2$ where $f_1$ and $f_2$ are bounded. Because $\phi_y\geq 0$,  we have that 
	
	\begin{multline}
		\frac{1}{\phi_{\varepsilon,y}(t,\nu_\varepsilon+y)^2} - \frac{1}{\phi_{\varepsilon,y}(t,\nu_\varepsilon-y)^2} = \frac{f_1(t,y) y^3 + f_2(t,y)(t_\varepsilon-t)^2}{\phi_{\varepsilon,y}(t,\nu_\varepsilon+y)\phi_{\varepsilon,y}(t,\nu_\varepsilon-y)^2} \\
		+\frac{f_1(t,y) y^3 + f_2(t,y)(t_\varepsilon-t)^2}{\phi_{\varepsilon,y}(t,\nu_\varepsilon+y)^2\phi_{\varepsilon,y}(t,\nu_\varepsilon-y)}.
	\end{multline}
	
	The two involved terms are similar. We consider only the first one. We obtain
	
	\begin{multline}
		|D_2| \leq \frac{C_\varepsilon}{(t_\varepsilon-t)^{1/2-2\lambda}} \int_{r=0}^{\kappa/2} \int_{z=0}^r \frac{1}{z^{1/2-2\lambda}} \\
		\cdot \frac{(r+z) (r-z)^4}{\left( (r+z)^2 + (t_\varepsilon-t) \right)^2\cdot \left( (r-z)^2 + (t_\varepsilon-t) \right) \cdot \left( (z-r)^2 + (t_\varepsilon-t) \right)^2}.
	\end{multline}
	
	Again, we split the domain of $z$ in two parts.
	
	\begin{multline}
		|D_2| \leq \frac{C_\varepsilon}{(t_\varepsilon-t)^{1/2-2\lambda}} \int_{r=0}^{\kappa/2} \int_{z=0}^{r/2} \frac{1}{z^{1/2-2\lambda}} \\
		\cdot \frac{(r+z) (r-z)^4}{\left( (r+z)^2 + (t_\varepsilon-t) \right)^2\cdot \left( (r-z)^2 + (t_\varepsilon-t) \right) \cdot \left( (z-r)^2 + (t_\varepsilon-t) \right)^2} \\
		+ \frac{C_\varepsilon}{(t_\varepsilon-t)^{1/2-2\lambda}} \int_{r=0}^{\kappa/2} \int_{z=r/2}^{r} \frac{1}{z^{1/2-2\lambda}} \\
		\cdot \frac{(r+z) (r-z)^4}{\left( (r+z)^2 + (t_\varepsilon-t) \right)^2\cdot \left( (r-z)^2 + (t_\varepsilon-t) \right) \cdot \left( (z-r)^2 + (t_\varepsilon-t) \right)^2}.
	\end{multline}
	
	Now, estimating each term separately, we obtain that $|D_2| \leq \frac{C_\varepsilon}{(t_\varepsilon-t)^{9/4-3\lambda}}$. Lastly, for $D_3$, we have to study the term $\chi'(\nu_\varepsilon-y) - \chi'(\nu_\varepsilon+y)$. We obtain
	
	\begin{equation}
		\left| \chi'(\nu_\varepsilon-y) - \chi'(\nu_\varepsilon+y) \right| \leq C_\varepsilon \cdot |y|.
	\end{equation}
	
	With similar computations, the additional $y$ leads to a gain of order $(t_\varepsilon-t)^{1/2}$. Lastly, for the terms $J_2$, $J_3$, $J_4$, the additional terms $x_1$, $y$, $y$ and $x_1 \cdot y$ converts into (respectively) a gain of order $(t_\varepsilon-t)^{1/2}$, $(t_\varepsilon-t)^{1/2}$, $(t_\varepsilon-t)^{1/2}$ and $(t_\varepsilon-t)^{1}$.
	
	Overall, we obtain that 
	
	\begin{equation}
		\int\int_{\beta\cup \delta} i(x_1,y) \geq \frac{C_\varepsilon}{\left( t_\varepsilon-t \right)^{11/4-3\lambda}}
	\end{equation}
	
	and hence 
	
	\begin{equation}\label{unresultat}
		(iv) \geq \frac{C_\varepsilon}{(t_\varepsilon-t)^{11/4-3\lambda}}.
	\end{equation}

	For $(i)$, using the lower bound on $\phi_\varepsilon$ provided by (\ref{properties3}), we obtain with the same method (we make the same change of variable and do not separate the domain as previously)
	
	\begin{equation}\label{i1lol}
		(i) \leq \int_x \int_y \frac{C_\varepsilon}{\left(t_\varepsilon-t\right)^{1/2-2\lambda}} \frac{M_\varepsilon |x|^{-1/4}}{x^2 + (t_\varepsilon-t)} \frac{M_\varepsilon |y|^{-1/4}}{y^2+(t_\varepsilon-t)} \leq \frac{C_\varepsilon}{(t_\varepsilon-t)^{14/8-2\lambda}}
	\end{equation}
	
	Also, for $(ii)$, we obtain
	
	\begin{multline}\label{i2lol}
		(ii) \leq \int_x \int_y \frac{C}{\left(t_\varepsilon-t\right)^{1/2-2\lambda}} \frac{M_\varepsilon |x|^{-1/4} }{x^2 + (t_\varepsilon-t)} \frac{M_\varepsilon |y|^{3/4} }{(y^2 + (t_\varepsilon-t))^2} \leq \frac{C_\varepsilon}{(t_\varepsilon-t)^{1/2+5/8+9/8-2\lambda}} \\
		\leq \frac{C_\varepsilon}{(t_\varepsilon-t)^{9/4-2\lambda}}
	\end{multline}
	
	The case of $(iii)$ is similar.
	Now, using (\ref{integrand3}), (\ref{unresultat}), (\ref{i1lol}) and (\ref{i2lol}), we obtain that 
	\begin{equation}
		I_\varepsilon^2(t) = (i) + (ii) + (iii) + (iv) \geq \frac{C_\varepsilon}{(t_\varepsilon-t)^{11/4-3\lambda}}
	\end{equation}
	
	Because of (\ref{imaj}), we get the desired result.
	
	\begin{equation}
		I_\varepsilon(t) = I_{\varepsilon}^1(t) + I_\varepsilon^2(t) + I_\varepsilon^3(t) \xrightarrow[t\to t_\varepsilon]{} \infty.
	\end{equation}
	
\end{proof}

\begin{remark}
	We also obtained a lower bound for the speed at which $||v||_{H^{7/4}} \rightarrow \infty$.
\end{remark}

\begin{theorem}
	$t_\varepsilon \rightarrow 0$ as $\varepsilon \rightarrow 0$.
\end{theorem}

\begin{proof}
	
	\begin{equation}\label{lifespan}
		\max \left\{ h_\varepsilon(y)=\frac{|\chi_\varepsilon'(y)|}{(1-\chi_\varepsilon(y))^2} \right\} \geq h_\varepsilon(\varepsilon) \geq |\ln(\varepsilon)|^\alpha \rightarrow \infty.
	\end{equation}
	
	By definition of $t_\varepsilon$, we thus have that $t_\varepsilon \rightarrow 0$ as $\varepsilon \rightarrow 0$.
	
\end{proof}

\section{A scaling argument to construct the solution}\label{chapterL}

In the following, we use a scaling argument to create a sequence of solutions with summable $H^{11/4}$ norms, and we choose a $\varepsilon$ parameter such that the lifespan goes to $0$.

Let $v(t,x)$ be a solution of (\ref{modelequation}). We define for $\omega, \gamma \in \mathbb{R}$,
\begin{equation}
	v_{\lambda} (t,x) = \lambda^\omega v(\lambda^\gamma t, \lambda^\gamma x).
\end{equation}

Now, 
\begin{equation}
	\begin{split}
		(\Box v_{\lambda})(t,x) &= \lambda^\omega \lambda^{2\gamma} \Box v(\lambda^\gamma t, \lambda^\gamma x) \\
		(Dv_{\lambda} D^2v_{\lambda})(t,x) &= \lambda^{2\omega} \lambda^{3\gamma} (Dv)(\lambda^\gamma t, \lambda^\gamma x) (D^2 v)(\lambda^\gamma t,\lambda^\gamma x);
	\end{split}
\end{equation}

so $v_\lambda$ is also a solution, provided that 
\begin{equation}\label{refi}
	\omega + \gamma = 0.
\end{equation}

Now, we will estimate the modification of the Sobolev norm due to the rescaling. Note that without the logarithmic modification, we have (with the change of variables in $x$ and $\xi$)

\begin{equation}
	\begin{split}
		||v_\lambda||_{\dot H^{11/4}} &= \lambda^\omega (\lambda^\gamma)^{11/4-2/2} ||v||_{\dot H^{11/4}}.
	\end{split}
\end{equation}

Now, we estimate $||v_\lambda||_{\dot H^{11/4}(\ln H)^{-\beta}}$. We have by definition

\begin{multline}
	||v_\lambda||_{H^{11/4}(\ln H)^{-\beta}}^2=  \int_{\xi \in \mathbb{R}^2} |\lambda|^{2\gamma} |\lambda|^{-4\gamma} \frac{|\lambda|^{2\omega} |\lambda^\gamma \xi|^{11/2}}{(1+|\ln(\lambda^\gamma |\xi|)|)^{2\beta}} \mathcal{F}(v)(\xi)^2 \\
	\leq (\lambda)^{2\omega} (\lambda^\gamma)^{11/2-2} \int_{\xi \in \mathbb{R}^2} \frac{|\xi|^{11/2}}{(1+|\ln(|\lambda^\gamma \xi|)|)^{2\beta}} \mathcal{F}(v)(\xi)^2,
\end{multline}

We have the following properties for $\lambda<1$,

\begin{equation}\label{deslemmes}
	\begin{split}
		&|\ln(|\lambda^\gamma \xi|)| = | \ln(\lambda^\gamma) + \ln(|\xi|)| \geq |\ln(\xi)| \text{ for } |\xi|\in (0,1], \\
		&\text{ for } |\xi| \in [1,\lambda^{-\gamma/2}],~ |\ln(\lambda^\gamma) + \ln(|\xi|)| \geq \frac{1}{2} |\ln(\lambda^\gamma)| \geq |\ln(|\xi|)|, \\
		&\text{ for } |\xi| \in [\lambda^{-2\gamma}, \infty),~ |\ln(\lambda^\gamma) + \ln(|\xi|)| \geq \frac{1}{2} |\ln(|\xi|)|.
	\end{split}
\end{equation}

This means that we can already establish 

\begin{multline}\label{lamblamb}
	\int_{|\xi|<\frac{1}{\lambda^{\gamma/2}} \bigcup |\xi|> \frac{1}{\lambda^{2\gamma}} } \lambda^{2\gamma} |\lambda|^{-4\gamma} \frac{|\lambda|^{2\omega} |\lambda^\gamma \xi|^{11/2}}{(1+|\ln(\lambda^\gamma |\xi|)|)^{2\beta}} \mathcal{F}(v)(\xi)^2 \\
	\leq 2 (\lambda)^{2\omega} (\lambda^\gamma)^{11/2-2} ||v||_{\dot H^{11/4} (\ln H)^{-\beta}}.
\end{multline}

For the last part, we have

\begin{multline}\label{lamblamb2}
	\int_{\lambda^{-\gamma/2} < |\xi|< \lambda^{-2\gamma}} \lambda^{2\gamma} \lambda^{-4\gamma} \frac{\lambda^{2\omega} |\lambda^\gamma \xi|^{11/2}}{(1+|\ln(\lambda^\gamma |\xi|)|)^{2\beta}} \mathcal{F}(v)(\xi)^2 \\
	\leq \int_{\lambda^{-\gamma/2} < |\xi|< \lambda^{-2\gamma}} \lambda^{2\gamma} \lambda^{-4\gamma} \lambda^{2\omega} |\lambda^\gamma \xi|^{11/2} \mathcal{F}(v)(\xi)^2 \\
	\leq \lambda^{-2\gamma} \lambda^{11/2 \gamma} \lambda^{2\omega} (1 + \left| \ln(\lambda^{2\gamma})\right|)^{2\beta} 
	\int_{\lambda^{-\gamma/2} < |\xi|< \lambda^{-2\gamma}} \frac{ | \xi|^{11/2}}{(1+ \left|\ln(|\xi|)\right|)^{2\beta}} \mathcal{F}(v)(\xi)^2  \\
	\leq \lambda^{2\omega} \lambda^{7/2 \gamma} (1+\left| \ln(\lambda^{2\gamma}) \right|)^{2\beta} ||v||_{\dot H^{11/4} (\ln H)^{-\beta}}.
\end{multline}

Overall, we obtain from (\ref{lamblamb}) and (\ref{lamblamb2}),

\begin{equation}
	||v_\lambda||_{\dot H^{11/4}(\ln H)^{-\beta}} \leq 2 \lambda^{\omega} |\lambda|^{\frac{7}{4} \gamma} (1+2 \left| \ln(\lambda^\gamma) \right|)^{\beta} ||v||_{\dot H^{11/4}(\ln H)^{-\beta}}.
\end{equation}

Now, we will choose the following values for the parameters and apply this result to the solution $u_\varepsilon$ defined in the previous chapter as the solution of the Cauchy problem (\ref{Cauchyjo}).

Now, applying this to the solution $u_\varepsilon$ defined in the previous chapter, define 
\begin{equation}
	\left\{
	\begin{split}
		\omega &= -1, \\
		\gamma &= 1, \\
		\lambda &= n^{-4}, \\
		\varepsilon &= \min \left(e^{-n^{5}},e^{-n^{\frac{6}{\alpha}}}\right),
	\end{split}
	\right.
\end{equation}

and $(u_n)$ the corresponding sequence of solutions.

First, the rescaling in time shortens the lifespan of $u_{\varepsilon}$ and we have that the lifespan of $u_{\varepsilon,\lambda}$, satisfies $t^{\lambda,\varepsilon} = \lambda^{-\gamma} t_\varepsilon$. Hence, we have by (\ref{lifespan})

\begin{equation}
	t_n \leq \frac{n^4}{|\ln(\varepsilon)|^\alpha} \leq \frac{1}{n} \rightarrow 0.
\end{equation}

Now, we can choose any parameter $\beta > 1/2$ in our construction. So we can assume for instance that we have $\beta< 3/4$.

\begin{equation}
	||u_{\lambda,\varepsilon}||_{H^{11/4} (\ln H)^{-\beta}} \leq 2 \frac{1}{n^3} (1+2 \ln(n^{-4}))^\beta ||u_{0,\varepsilon}||_{H^{11/4} (\ln H)^{-\beta}},
\end{equation}

which is in $l^1(\mathbb{Z})$.

Now, define $\tilde{u}_n$ as a translation of $u_{n}$ in $x_1$ such that $\text{Supp}(\tilde u_i) \cap \text{Supp}(\tilde u_j) = \emptyset$ for $i\neq j$ (and the domain of dependence do not intersect either), and finally define 

\begin{equation}\label{formuleL}
	\Xi(t,x) = \sum_{n=0}^{\infty} \tilde{u}_n(t,x).
\end{equation}

Each function $\tilde{u}_n$ is a translation of a function of which the support is included in $x_1 \in [\varepsilon/2, \frac{2}{|\ln(\varepsilon)|^{\alpha/2}}] \subseteq [\frac{1}{2} e^{-n^5}, \frac{2}{n^3}]$. Because the sequence $\frac{2}{n^3} - \frac{1}{2}e^{-n^5}$ is in $l^1(\mathbb{N})$, we can find a sequence of translations such that the projection of the support of $\Xi$ on the $x_1$ axis is bounded. Now, the width of the domain (the projection on the $x_2$ axis) is bounded by a constant because it is bounded by $\sqrt{\frac{1}{|\ln(\varepsilon_n)|^{\alpha/2}}} \cdot \left|\ln\left(\frac{1}{|\ln(\varepsilon_n)^{\alpha/2}}\right)\right|^{-\delta} \rightarrow 0$ as $n\rightarrow \infty.$ This means that $\Xi$ has a compact support. In virtue of lemma \ref{celuila}, we hence obtain

\begin{theorem}
	$\Xi$ satisfies
	\begin{equation}
		\begin{split}
			||\Xi_{|t=0}||_{H_{x_1}^{11/4}} &< \infty \\
			||\frac{\partial \Xi_{|t=0}}{\partial t}||_{H_{x_1}^{7/4}} &< \infty \\
			\forall t>0, \hspace{0.2cm} ||\frac{\partial}{\partial x} (D\Xi) (t,\cdot)||_{H_{x_1}^{3/4}} &= \infty 
		\end{split}
	\end{equation}
\end{theorem}

\begin{proof}
	It directly follows from the fact that $t_n \rightarrow 0$ as $n\rightarrow \infty$. 
\end{proof}

\begin{remark}
	The function $\Xi$ that we created to show the ill-posedness of the equation, is also of compact support.
\end{remark}

\nocite{*}
\bibliographystyle{alpha}
\bibliography{bibliography}

\appendix 

\section{Detailed proof of the first case of prop 4.1}\label{appendixA}

	Let us now prove \ref{35lol}. We assume that $|v| < 1/100$, and that $t\leq C \frac{1}{|\ln(\varepsilon)|^\alpha}$ (This will be achieved whenever $\varepsilon$ is small enough, which means we will only consider times $t \leq  C \frac{1}{|\ln(\varepsilon)|^\alpha}$ very small, see part 3 for further details.)
	
	We will distinguish three cases, first, we consider curves whose starting point has an abscissa strictly bigger than $x_0 = \frac{\varepsilon}{4}$.
	
	It follows from definition \ref{dependence} that $(t',x') \in \Omega$ if and only if $(t',x_1') \in \Omega^1$ and all Lipschitz continuous curves from $(t',x')$ that satisfy (\ref{condlight}) intersect the hyperplane $t=0$ in the set $\{ x | |x_2| \leq \frac{\sqrt{2x_1}}{|\ln(x_1)|^\delta} \}$. 
	
	Now, let $(t(s),x_1(s),x_2(s))$ be a Lipschitz continuous curve parameterized so that $t(s) + x_1(s) = s$. Note $q(s) = x_1(s) - t(s)$. Note that (\ref{condlight}) is equivalent to (using the fact that $\frac{dt(s)}{ds}+\frac{dx_1(s)}{ds}=1$),

	\begin{equation}
		R(s) \leq v(t(s),x_1(s)) - \frac{d q(s)}{ds},
	\end{equation}
	where $R(s) = \left(\frac{dx_2(s)}{ds}\right)^2$.
	
	Now, using this set of new variables $s= x_1 +t$, $q = x_1-t$ and $U(s,q) = u((s-q)/2,(s+q)/2)$, (\ref{dimone}) becomes
	
	\begin{equation}
		\left\{
		\begin{split}
			(\partial_s + V(s,q) \partial_q) \partial_q U(s,q) &= 0, \hspace{0.2cm} V(s,q) = 2 \partial_q U(s,q),\\
			U(y,y) &= 0, \hspace{0.2cm} U_q(y,y) = \frac{1}{2}\chi(y)
		\end{split}
		\right.
	\end{equation}
	
	The characteristics are given by $s = constant$ and $q =h(s,y)$ with 
	\begin{equation}
		\frac{d}{ds}h(s,y) = V(s,h(s,y)), \hspace{0.2cm} h(y,y)= y.
	\end{equation}
	Thus, $s \mapsto V(s,h(s,y))$ is constant on the curve and is equal to $\chi(y)$. These are the key ingredients to make this proof. 
	
	Now assume that the curve is such that $q(a) = a$ and $q(b) = h(b,y)$. With $r(s) = \sqrt{x_2^2}$, assume in addition that $r(a)^2 = \frac{2a}{|\ln(a)|^{2 \delta}}$. Note that because $q(a) = a$ is equivalent to $t(a) = 0$, $r(a)^2 = \frac{2a}{|\ln(a)|^{2 \delta}}$ is simply that the point $(t(a),x_1(a),x_2(a)) = (0,a, \frac{\sqrt{2 a}}{|\ln(a)|^{\delta}})$ is on the edge of the domain $\Omega_0$.
	Also, the condition $q(b) = h(b,y)$ is equivalent to $x_1(b) = \phi(t(b),y)$.
	
	Now, 
	\begin{equation}
		\begin{split}
			|r(b) - r(a)| \leq \int_{a}^b \sqrt{R(s)} ds \leq \sqrt{b-a} \sqrt{\int_a^b R(s)ds}  \\
			\leq \sqrt{b-a} \sqrt{q(a) - h(b,y)} \hspace{0.2cm}\text{ (Using $V_q \leq 0$) } 
		\end{split}
	\end{equation}
	
	Now, 
	\begin{multline}
		q(a) - h(b,y) = a - y + h(y,y) - h(b,y) = a- y + \int_y^a V(s,h(s,y)) ds\\
		\leq \frac{101}{100}(a-y) \leq \frac{9}{8}a
	\end{multline}
	
	Also, 
	\begin{multline}
		b-y =2t + h(b,y) - h(y,y) = 2t + \int_y^b V(s,h(s,y)) ds \leq 2t + \frac{b-y}{100}
	\end{multline}
	Hence $\sqrt{b-a} \leq \sqrt{b-y} \leq \sqrt{2} \cdot \sqrt{t}$, because $a \geq y$. (It comes from the expression of $h$ using $\phi$.)

	Now, since $a = \frac{1}{2}r(a)^2 |\ln(a)|^{2\delta}$, we obtain that $a \geq \varepsilon/4$ implies $|\ln(a)|^\delta \leq C |\ln(\varepsilon)|^\delta$. Hence, using condition (\ref{condalphadelta}), we obtain
	
	\begin{equation}
		|r(b)-r(a)| \leq  C \sqrt{t} \cdot \sqrt{a} \leq C r(a) \sqrt{t} \cdot |\ln(\varepsilon)|^{\delta} \leq r(a) |\ln(\varepsilon)|^{\delta-\alpha/2} \leq \frac{1}{2} r(a) , 
	\end{equation}
	
	\begin{equation}
		r(b) \geq r(a) - |r(a)-r(b)| \geq r(a) - 2/3 r(a)  \geq 1/2 r(a).
	\end{equation}
	
	Meaning that any curve starting with an abscissa bigger than $x_0$ does not reach the inside of a ball centered in $(t=t, x_1 = \phi(t,\nu_\varepsilon), x_2= 0)$ of a certain radius $\delta_1$, for a fixed $\nu_\varepsilon$ only depending on $\varepsilon$. However, the radius $\delta_1$ may depend on $\varepsilon$.

\end{document}